\documentclass{scrartcl}

\usepackage{graphicx}
\usepackage[table]{xcolor}

\usepackage{booktabs,units,xfrac}
\usepackage{placeins}
\usepackage{pgfplots}
\usepackage{siunitx}
\pgfplotsset{compat=1.9}

\usetikzlibrary{
  pgfplots.groupplots,
  matrix
}
\usepackage{amsmath,amssymb,dsfont,graphicx}

\usepackage[table]{xcolor}

\usepackage{algpseudocode}
\usepackage{algorithm}
\usepackage{float}
\floatstyle{plaintop}

\usepackage{tikz}
\usepackage{pgfplots}
\newlength \figureheight
\newlength \figurewidth

\newcommand{\R}{\mathds{R}}

\newcommand{\tr}{\operatorname{tr}}

\newcommand{\vt}{\mathbf{v}}
\newcommand{\pt}{\mathbf{p}}

\newcommand{\ut}{\mathbf{u}}
\newcommand{\zt}{\mathbf{z}}
\newcommand{\bt}{\mathbf{b}}
\newcommand{\Ft}{\mathbf{F}}
\newcommand{\Et}{\mathbf{E}}

\newcommand{\deltat}{\mathbf{\delta}}

\newcommand{\sigmat}{\boldsymbol{\sigma}}
\newcommand{\Sigmat}{\boldsymbol{\Sigma}}

\makeatletter
\newcommand{\opnorm}{\@ifstar\@opnorms\@opnorm}
\newcommand{\@opnorms}[1]{%
  \left|\mkern-1.5mu\left|\mkern-1.5mu\left|
   #1
  \right|\mkern-1.5mu\right|\mkern-1.5mu\right|
}
\newcommand{\@opnorm}[2][]{%
  \mathopen{#1|\mkern-1.5mu#1|\mkern-1.5mu#1|}
  #2
  \mathclose{#1|\mkern-1.5mu#1|\mkern-1.5mu#1|}
}
\makeatother

\newcommand{\SO}{\mathcal{S}}
\newcommand{\IN}{\mathcal{I}}

\newcommand{\FL}{\mathcal{F}}

\newcommand{\norm}[1]{\lVert #1 \rVert}

\newcommand{\SCC}{\cellcolor{gray!20}}

\definecolor{mycolor1}{rgb}{0.00000,0.44700,0.74100}%
\definecolor{mycolor2}{rgb}{0.85000,0.32500,0.09800}%
\definecolor{mycolor3}{rgb}{0.92900,0.69400,0.12500}%
\definecolor{mycolor4}{RGB}{152,198,234}%
\definecolor{mycolor5}{RGB}{1,1,1}%

\title{A Newton multigrid framework for optimal control of
  fluid-structure interactions
}

\author{L. Failer\thanks{%
    Technische Universit\"at M\"unchen, 
    85748 Garching bei M\"unchen, Germany, 
    \texttt{lukas.failer@ma.tum.de}}
  \and
  T. Richter\thanks{%
    Otto-von-Guericke Universit\"at Magdeburg, 
    39104 Magdeburg, Germany, 
    \texttt{thomas.richter@ovgu.de}, 
    and Interdisciplinary Center for Scientific Computing, 
    Heidelberg University, 69120 Heidelberg, Germany}}

\begin{document}

\maketitle

\begin{abstract}
  In this paper we consider optimal control of nonlinear time-dependent
fluid structure interactions. To determine a time-dependent control
variable a BFGS algorithm is used, whereby gradient information is
computed via a dual problem. To solve the resulting ill conditioned
linear problems occurring in every time step of state and dual
equation, we develop a highly efficient monolithic solver that is
based on an approximated Newton scheme for the primal equation and a
preconditioned Richardson iteration for the dual problem. The
performance of the presented algorithms is tested for one 2d and one 3d
example numerically.  

Keywords:
  fluid-structure interactions; finite elements; 
  multigrid; optimal control; parameter estimation
\end{abstract}

\section{Introduction}

  Fluid-structure interactions are part of various applications
  ranging from classical engineering problems like aeroelasticity or
  naval design to medical applications, e.g. the flow of blood in the
  heart or in blood vessels. More and more of these applications are regarded recently in combination with optimal control, shape-optimization, and parameter estimation. Especially in hemodynamical applications --- in order to get a deeper understanding of the development of vascular diseases --- patient specific properties have to be incorporated into the models. For example, in~\cite{BertoglioBarberGaddum,BertoglioChapelleFernandezGerbeauMoireau,parameterhemodynamics,DEliaMirabellaPasserini,LassilaManzoniQuarteroniRozza,PantFabregesGerbeau,MoireauBertoglioXiao} patient specific boundary conditions and vessel material parameters are determined to simulate arterial blood flow. Similar approaches using gradient information have been proposed in ~\cite{DEliaMirabellaPasserini,BertagnaDEliaPerego,Cardiovascular} 
  to estimate Young's modulus of an artery. 
 
  As computer tomography (CT) and magnetic resonance imaging (MRI)
  evolve rapidly,  already very accurate measurements of the movement
  of the vessel wall are possible nowadays and even averaged flow
  profiles in blood vessels can be provided, see
  \cite{AsnerHadjicharalambous,BertoglioBarberGaddum,LamataPitcherNordsletten}. To
  incorporate the data in the vascular models, it is necessary to
  improve the available parameter estimation and optimal control
  algorithms for fluid-structure interaction applications, in
  particular since only few approaches in the literature take the
  sensitivity information of the full time-dependent nonlinear system
  into account. For example
  in~\cite{DegrooteHojjatStavropoulouWuechnerBletzinger,MartinClementDecoeneGerbeau},
  adjoint equations are derived for one-dimensional fluid-structure
  interaction configurations and in~\cite{RichterWickOpt} for a
  stationary  fluid-structure interaction problem. In contrast, the
  authors
  of~\cite{PantFabregesGerbeau,parameterhemodynamics,MoireauBertoglioXiao,BertoglioBarberGaddum}
  use a sequential reduced Kalman filter. Peregio, Veneziani, and
  Vergara~\cite{Cardiovascular} compute sensitivity information to
  estimate the wall stiffness. To reduce the computational time, they
  solve in every time-point an optimal control problem. As the mesh
  motion is discretized via an explicit time-stepping scheme, no
  sensitivity information  of the mesh motion equation has to be
  computed. Similar to the articles
  \cite{parameterhemodynamics,BertoglioBarberGaddum,MoireauBertoglioXiao},
  the estimated parameters are updated in every time step and the forward simulation only runs once. 
  
  In this paper we are going to compute gradient information for the
  full time-dependent nonlinear system for 3d applications. Thereby,
  the optimization algorithm takes the intrinsic property of
  fluid-structure interaction, transport over time, into account. In
  addition the here presented approach enables to regard tracking type
  functionals with observation at a singular time-point or on a
  specific time-interval. Furthermore a time-dependent parameter can
  be reconstructed. This would not be possible, if we would use a
  Kalman filter or would solve an optimization problem in every time step as in the literature cited above. The dual problem to compute sensitivities can be derived as in \cite{Failer2017} or as in \cite{FailerWick}, where sensitivity information was used for a dual-weighted residual error estimator. 
 
  For various applications
  and a general overview on 
  modeling and discretization techniques for fluid-structure
  interactions we refer to the
  literature~\cite{BazilevsTakizawaTezduyar2013,Richter2017}. Mathematically, two challenges come
  together in fluid-structure interactions: First, fluid-structure interactions are free boundary value problems.
  The governing domains for the fluid - we will
  consider the incompressible Navier-Stokes equations - and the solid
  - we consider hyperelastic materials like the St. Venant Kirchhoff
  model - move and the motion is determined by the coupled dynamics,
  i.e. it is not known a priori. This geometric
  problem is treated by mapping onto a fixed domain~\cite{Donea1982} such that movement of the boundary is incorporated into the equation and we can derive the dual problem on the fixed reference domain.
  Second, the two problems that are
  coupled are of different type, the parabolic Navier-Stokes equation
  and the hyperbolic solid problem. On the common and moving
  interface, both systems are coupled by different conditions. This
  coupling gives rise to stability problems that can call for small
  time steps or many subiterations. Most prominently this problem shows itself in the so called
  added mass effect~\cite{CausinGerbeauNobile2005}. The added mass effect is of
  particular relevance in hemodynamical applications, that are focus of
  this work~\cite{HronTurekMadlikEtAl2010}, and calls for monolithic
  formulations and strongly coupled discretizations and solution techniques. 
  This property is transmitted to the dual problem, such that we have to derive strongly coupled solution techniques for the dual problem. 
  To compute dual information, we extend the Newton solver proposed
  for time-dependent fluid-structure interactions in
  \cite{FailerRichter2019} to the dual problem. Thereby, iterative
  solvers, preconditioned with geometric multigrid, can solve the
  resulting linear problems in every state and dual time step very robust and efficiently.

In Section \ref{sec:eq} we present the optimal control problem, which is discretized in Section \ref{sec:disc} in space and time. For the discretized system we derive optimality conditions. In Section \ref{sec:solver} we discuss modifications for the Newton scheme presented first in \cite{FailerRichter2019} and extend the approach to the dual problem. Finally we test the proposed algorithm in Section \ref{sec:num} numerically to analyze the behavior of the Newton scheme. In addition we take a closer look on the convergence behavior of the iterative solvers.

\section{Governing equations}\label{sec:eq}

Here, we present the optimal control problem of a tracking type functional subject to fluid structure
interactions. We use a monolithic formulation for the fluid-structure interaction model coupling the incompressible Navier-Stokes equations and
an hyperelastic solid, based on the St. Venant Kirchhoff material. For
details we refer to~\cite{Richter2017}. The here presented
optimization approach can be directly extended to specific material
laws used in hemodynamics. As control variable we chose exemplarily the mean pressure over time at the outflow boundary. In the following we restrict us to the control space $Q=L^2(I)$, but the here presented optimization algorithm can as well be applied to determine material parameters (e.g.~$Q=\R^n$) or space-distributed parameters (e.g.~$Q=L^2(\Omega)$) entering the fluid- or solid-problem. 

On the $d$-dimensional domain, partitioned in reference configuration
$\Omega = \FL\cup\IN\cup\SO$, where $\FL$ is the fluid domain, $\SO$
the solid domain and $\IN$ the fluid structure interface,  we denote
by $\vt$ the velocity field, split into fluid  velocity
$\vt_f:=\vt|_{\FL}$ and solid velocity $\vt_s:=\vt|_{\SO}$, and 
by $\ut$ the deformation field, again with $\ut_s:=\ut|_{\SO}$ and
$\ut_f:=\ut|_{\FL}$. The boundary of the fluid domain
$\Gamma_f:=\partial\FL\setminus\IN$ is split into inflow boundary
$\Gamma_f^{in}$ and wall boundary $\Gamma_f^{wall}$, where we usually
assume Dirichlet conditions,
$\Gamma_f^D:=\Gamma_f^{in}\cup\Gamma_f^{wall}$,  and a possible outflow boundary
$\Gamma_f^{out}$, where we enforce the do-nothing outflow
condition~\cite{HeywoodRannacherTurek1992}, and the control boundary $\Gamma_q$. The solid boundary
$\Gamma_s=\partial\SO\setminus\IN$ is split into Dirichlet part
$\Gamma_s^D$ and a Neumann part $\Gamma_s^N$.

We formulate the coupled fluid-structure interaction problem in a
strictly monolithic scheme by mapping the moving fluid domain onto the
reference state via the ALE map $T_f(t):\FL\to\FL(t)$,
constructed by a fluid domain deformation $T_f(t)=\operatorname{id} +
\ut_f(t)$. In the solid domain, this map
$T_s(t)=\operatorname{id}+\ut_s(t)$ denotes the Lagrange-Euler
mapping and as the deformation field $\ut$ will be defined globally on
$\Omega$ we simply use the notation $T(t)=\operatorname{id}+\ut(t)$
with the deformation gradient $\Ft:=\nabla T$ and its determinant
$J:=\operatorname{det}(\Ft)$. 

For given desired states $\tilde \vt(t) \in L^2(\FL)$ or $\tilde \ut(t) \in L^2(\SO)$, we find the global (in fluid and solid domain)
velocity and deformation fields
\begin{align*}
\vt(t)&\in \vt^D(t)+H^1_0(\Omega;\Gamma_f^D\cup\Gamma_s^D)^d  \text{ and } 
\ut(t)\in
\ut^D(t)+H^1_0(\Omega;(\partial\FL\setminus\IN)\cup\Gamma_s^D)^d,
\end{align*}
the  pressure $p\in L^2(\FL)$ and the control parameter $q\in Q$ 
satisfying the initial condition $\vt(0)=\vt_0$ and $\ut(0)=\ut_0$, 
as solution to 
\begin{equation}\label{aletime:1}
\begin{aligned}
\min_{q \in Q} J(q,\vt,\ut)= \frac{1}{2}\int_{I} \norm{\vt-\tilde \vt}^2_{\FL} \text{ d}t+ \frac{1}{2}\int_{I} \norm{\ut-\tilde \ut}^2_{\SO}\text{ d}t +\frac{\alpha}{2}\norm{q}^2_Q
\end{aligned}
\end{equation}
and subject to
\begin{equation}\label{aletime:1.1}
  \begin{aligned}
    \big( J(\partial_t  \vt +  (\Ft^{-1}( \vt-\partial_t   
    \ut)\cdot\nabla) \vt,\phi\big)_{\FL}
    + \big( J\sigmat_f
    \Ft^{-T},\nabla\phi\big)_{\FL} \\
   +(\rho_s^0\partial_t \vt,\phi)_{\SO}
    +( \Ft\Sigmat_s,\nabla\phi)_{\SO}
    &=
    (q,\phi)_{\Gamma_q} \\
    \big(J\Ft^{-1}:\nabla\vt^T,\xi\big)_{\FL} &= 0\\
    (\partial_t \ut- \vt,\psi_s)_{\SO}&=0\\
    (\nabla \ut,\nabla\psi_f)_{\FL}&=0,
  \end{aligned}
\end{equation}
where the test functions are given in 
\[
\phi\in H^1_0(\Omega;\Gamma_f^D\cup\Gamma_s^D)^d,\quad
\xi\in L^2(\FL),\quad
\psi_f\in H^1_0(\FL)^d,\quad
\psi_s\in L^2(\SO)^d.
\]
By $\rho_s^0$ we denote the solid's density, by $\ut^D(t)\in
H^1(\Omega)^d$ and $\vt^D(t)\in H^1(\Omega)^d$  extensions of the
Dirichlet data into the domain.  
The Cauchy stress
tensor of the Navier-Stokes equations in ALE coordinates is given by
\[
\sigmat_f(\vt,p) = -p_f I + \rho_f\nu_f (\nabla\vt\Ft^{-1} + \Ft^{-T}\nabla\vt^T)
\]
with the kinematic viscosity $\nu_f$ and the density $\rho_f$. In the
solid we consider the St. Venant Kirchhoff 
material with the second Piola Kirchhoff tensor $\Sigmat_s$ based on the Green
Lagrange strain tensor $\Et_s$ 
\[
\Sigmat_s(\ut) = 2\mu_s \Et_s + \lambda_s\tr(\Et_s)I,\quad
\Et_s:=\frac{1}{2}(\Ft^T\Ft-I)
\]
and with the shear modulus $\mu_s$ and the Lam\'e coefficient
$\lambda_s$. In~(\ref{aletime:1.1}) we construct the ALE extension
$\ut_f=\ut|_\FL$ by a simple harmonic extension. A detailed discussion
and further literature on the construction of this extension is found
in~\cite{YirgitSchaeferHeck2008,Richter2017}. For shorter notation, we denote by $U:=(\vt,\ut,p_f)\in X$ the solution variable and with $X$ the corresponding ansatz space 
and by $\Phi:=(\phi,\psi_f,\psi_s,\xi)\in Y$ the test functions and the corresponding test space.

For a control $q\in L^2(I)$ and the here given tracking-type functional constrained by linear-fluid structure interaction, we were able to proof in \cite{FailerMeidnerVexler} existence of a unique solution and $H^1(I)$ regularity of the optimal control. In addition an optimality system could be rigorously derived. Due to the missing regularity results for the here regarded nonlinear control to state mapping, no further theoretical conclusions are possible here.

\section{Discretization}\label{sec:disc}

In the following we give a description of the discretization of the
fluid-structure interaction system~(\ref{aletime:1.1}) in space and in 
time. While there exist many variants and different realizations, our
choice of methods is based on the following principles
\begin{itemize}
\item Since the fsi system is a constraint in the optimization process
  we base the discretization on Galerkin methods in space and
  time. This helps us to derive the discrete optimality system. As far
  as possible (up to quadrature error) we aim at permutability of
  discretization and optimization.
\item Aiming at three dimensional problems we consider methods of
  reasonable approximation error at feasible costs. In space we will
  use second order finite elements and in time a second order time
  stepping scheme. This approach is similar
  to~\cite{HronTurekMadlikEtAl2010} or our previous work documented in
  \cite{Richter2017}. 
\item Since the key component of the linear solver is a geometric
  multigrid method with Vanka type blocking in the smoother we choose
  equal-order finite elements for all unknowns, pressure, velocity and
  deformation adding stabilization terms for the inf-sup
  condition. This setup allows for efficient linear algebra and  local
  blocking of the unknowns that is in favor of strong local 
  couplings taking care of all nonlinearities, see
  also~\cite{BraackRichter2006a} for a detailed description of the
  realization in the context of reactive flows. 
\item The temporal dynamics of fluid-structure interactions is
  governed by the parabolic/hyperbolic character of the
  coupling. In particular long term simulations give rise to stability
  problems. The Crank-Nicolson shows stability problems such that
  variants will be considered, see~\cite{RichterWick2015_time}. 
\end{itemize}

\subsection{Temporal discretization}

In~\cite{RichterWick2015_time} and~\cite[Section 4.1]{Richter2017}
many aspects of time discretization of monolithic fluid-structure
interactions are discussed. It turns out that the standard
Crank-Nicolson scheme is not sufficiently stable for long time
simulations. Suitable variants are the fractional step theta
method or shifted versions of the Crank-Nicolson scheme which
we refer to as theta time stepping methods. Applied to the ode
$u'=f(t,u(t))$ they take the form
\[
  u_n-u_{n-1} = k_n \theta f(t_n,u_n) + k_n(1-\theta) f(t_{n-1},u_{n-1}), 
\]
if $0=t_0<t_1<\cdots<t_N=T$ are the discrete time steps with step size
$k_n=t_n-t_{n-1}$. The choice $\theta=\frac{1}{2}+{\cal O}(k)$ gives
second order convergence and sufficient
stability~\cite{RichterWick2015_time}. Alternative approaches are the
fractional step theta scheme that consists of three sub steps with
specific choices for $\theta$ and the step size or the enrichment of
the Crank-Nicolson scheme with occasional Euler steps,
see~\cite{Rannacher1984}. 

In the context of optimization problems we aim at permutability of
optimization and discretization such that Galerkin approaches are of a
favor. In~\cite{MeidnerRichter1,MeidnerRichter2} we have
demonstrated an interpretation of the general theta scheme and the
fractional step theta scheme as Galerkin method with adapted function
spaces: the solution is found in the space of continuous and piecewise
(on $I_n=(t_{n-1},t_n)$) linear functions, the test-space is a space
rotated constant functions with jumps at the discrete time steps
$t_n$, namely
\[
  \psi^\theta\big|_{I_n}(t) = 1+ \frac{(6\theta-3)(2t-t_{n-1}-t_n)}{k_n}. 
\]
The theta scheme is recovered exactly for linear problems and
approximated by a suitable quadrature rule for nonlinear problems.

In case of fluid-structure interactions the domain motion term
$(J\Ft^{-1}\partial_t\ut\cdot\nabla\vt,\phi)$ takes a special role
since it couples temporal and spatial differential
operators. In~\cite{RichterWick2015_time} various discretizations are
analyzed and all found to give results in close agreement. 

Here, we consider the Galerkin variant of the theta scheme and we
approximate all temporal integrals by the quadrature rule
(see~\cite{MeidnerRichter2}) 
\[
\int_{t_{n-1}}^{t_n} f(t)\psi^\theta(t)\,\text{d}t
=k_n\theta f(t_n) + k_n(1-\theta)f(t_{n-1}) + {\cal O}\big(k_n^2
\| f\|_{W^{2,1}([t_{n-1},t_n])}\big). 
\]
The resulting discrete scheme is - up to quadrature error - the
standard theta time stepping scheme, which we use in our
implementation for reasons of efficiency. 

For the following we denote by $U_n\approx U(t_n)$ the approximation
at time $t_n$. Further we introduce
\begin{equation}\label{aletime:1.5}
  \begin{aligned}
    &A_F(U,\phi) := \big(J(\Ft^{-1}\vt\cdot\nabla)\vt,\phi\big)_\FL
    + \big(\rho_f\nu_f
    J(\nabla\vt\Ft^{-1}+\Ft^{-T}\nabla\vt^T)\Ft^{-T},\nabla\phi\big)_\FL \\
    &A_S(U,\phi) :=  \big(\Ft\Sigmat_s,\nabla\phi\big)_{\SO} ,\quad
    A_{ALE}(U,\psi_f) := \big(\nabla\ut,\nabla\psi_f\big)_{\FL}    
    \\
    &A_{p}(U,\phi) :=  \big(Jp\Ft^{-1},\nabla\phi\big)_{\FL} ,\quad
    A_{div}(U,\xi) := \big( J\Ft^{-1}:\nabla\vt^T,\xi\big)_{\FL}\\
    &F_{TR}(U_n,U_{n-1},\phi):=\big(( \bar J_n\bar\Ft^{-1}
    (\ut_n-\ut_{n-1})\cdot\nabla)\bar\vt_n,\phi\big)_{\FL}
    ,
  \end{aligned}
\end{equation}
and the  step $t_{n-1}\mapsto t_n$ is given as 
\begin{equation}\label{aletime:2}
  \begin{aligned}
    \big(\bar J_n (\vt_n-\vt_{n-1}),\phi\big)_{\FL}
     -F_{TR}(U_n,U_{n-1},\phi)
    + k A_p(U_n,\phi) +k\theta A_F(U_n,\phi)  \\
    +\big(\rho^0_s (\vt_n-\vt_{n-1}),\phi\big)_{\SO}+k\theta A_S(U_n,\phi)\\
    = - k(1-\theta) A_F(U_{n-1},\phi)- k(1-\theta) A_S(U_{n-1},\phi) +k(q_n,\phi)_{\Gamma_q}  \\
    k A_{div}(U_n,\xi) =0\\
    k A_{ALE}(U_n,\psi_f) =0\\
    \big(\ut_n,\psi_s\big)_{\SO}
    -k\theta\big(\vt_n,\psi_s\big)_\SO =
    \big(\ut_{n-1},\psi_s\big) + 
    k(1-\theta)\big(\vt_{n-1},\psi_s\big)_\SO,\\
  \end{aligned}
\end{equation}
with $\bar J_n = \sfrac{1}{2}(J_{n-1}+J_n)$ and $\bar\Ft_n =
\sfrac{1}{2}(\Ft_{n-1} + \Ft_n)$. The divergence equation $A_{div}$
and the pressure coupling $A_p$ are fully implicit, which can be
considered as a post processing step, see~\cite{MeidnerRichter1}.

If the optimality system is first derived and then discretized using
the Petrov-Galerkin discretization, we could observe that the
control variable has to be in the theta dependent test space of the
adjoint variable. As this space is very difficult to interpret  
the control variable $q\in Q$ is approximated by piece-wise constant
functions $q_n$ on every time-interval in the following. An
alternative interpretation is to actually use the theta dependent test
space for the adjoint variable but to approximate these integrals with
the midpoint rule giving
\[
\int_{t_{n-1}}^{t_n} f(t)\psi^\theta(t)\,\text{d}t
=k_n f(t_{n-\frac{1}{2}}) + {\cal O}\left(k_n\big|2\theta-1\big|
\|f\|_{W^{2,1}([t_{n-1},t_n])}\right). 
\]
Given the choice $\theta=\sfrac{1}{2}+{\cal O}(k_n)$ this gives
correct 
second order convergence. 
Numerical studies comparing both approaches did not result in a different
behavior of the optimization algorithm.  

\subsection{Finite elements}

Spatial discretization of the primal and adjoint problem is by means
of quadratic finite elements in all variables on a quadrilateral and
hexahedral meshes. The interface $\IN$ is resolved by the mesh such
that no additional approximation error appears. To cope with the
saddle point structure of the flow problem we use the local
projection method for
stabilization~\cite{BeckerBraack2001,Frei2016,Molnar2015,Richter2017}. In
the context of optimization problems this scheme has the advantage
that stabilization and optimization  commute, see~\cite{Braack2009}.
Further details on this and comparable approaches are found in the
literature~\cite{HronTurekMadlikEtAl2010,RichterWick2010,Richter2017}. 

The use of equal order finite elements in all variables has the
advantage that one set of scalar test functions
$\{\phi_h^{(1)},\dots,\phi_h^{(N)}\}$ can be chosen for all
variables. The discrete solution $U_h$ can then be written as
\[
  U_h(x) = \sum_{i=1}^N \mathbf{U}_i \phi_h^{(i)}(x)
\]
with coefficient vectors $\mathbf{U}_i=(p_i,\vt_i,\ut_i)\in
\R^{2d+1}$ and scalar test functions $\phi_h^{(i)}$. Likewise, the
resulting matrix entries 
$A_{ij}=A'(U_h)(W_h^{(j)},\Phi_h^{(i)})$ are small but dense local
matrices of size $(2d+1)\times (2d+1)$. All linear algebra routines
act on these blocks, e.g. inversion of a matrix entry corresponds to
the inversion of these blocks $A_{ij}^{-1}$, which results in a better
cache efficiency and reduced effort for indirect indexing of matrix
and vector entries. The effect of this approach is described
in~\cite{BraackRichter2006a}. 

\subsection{Optimality system and adjoint equation}\label{Sec::OptSystem}
As gradient based algorithms for parameter estimation are not very common in the hemodynamics community, we shortly derive the Karush-Kuhn-Tucker system and show how gradient information can thereby be extracted. 
To derive the Karush-Kuhn-Tucker system, we define Lagrange multipliers $Z_n=(\zt_n^p,\zt_n^v,\zt_n^{uf},\zt_n^{us})\in Y_h$ in every time step $n=0,...,N$ and get the discrete Lagrangian $L:\left(\R^N,(X_h)^{N+1},(Y_h)^{N+1}\right)\longmapsto \R$:
\begin{equation}\label{lagrangian}
  \begin{aligned}
    &L((q_n)_{n=1}^N,(U_n)_{n=0}^N,(Z_n)_{n=0}^N):=\\
    &\sum_{n=1}^{N-1} \Big\{\frac{1}{2} k \norm{\vt_n-\tilde \vt(t_n)}^2_{\FL} +\frac{1}{2} k \norm{\ut_n-\tilde \ut(t_n)}^2_{\SO}+ \frac{\alpha}{2} k q_n^2 \Big\}\\ 
& +\frac{1}{4} k \norm{\vt_0-\tilde \vt(t_0)}^2_{\FL} +\frac{1}{4} k \norm{\ut_0-\tilde \ut(t_0)}^2_{\SO}\\
& +\frac{1}{4} k \norm{\vt_N-\tilde \vt(t_N)}^2_{\FL} +\frac{1}{4} k \norm{\ut_N-\tilde \ut(t_N)}^2_{\SO}+\frac{\alpha}{2} k q_N^2\\
    &- \sum_{n=1}^{N} \Big\{
    \big(\rho^0_s (\vt_n-\vt_{n-1}),\zt^v_n \big)_{\SO}+k\theta A_S(U_n,\zt^v_n)+ k(1-\theta) A_S(U_{n-1},\zt^v_n) \\
    &\qquad+ \big(\ut_n,\zt^{us}_n\big)_{\SO} -\big(\ut_{n-1},\zt^{us}_n\big)-k\theta\big(\vt_n,\zt^{us}_n\big)_\SO -k(1-\theta)\big(\vt_{n-1},\zt^{us}_n\big)_\SO\\
    &\qquad+(\bar J_n (\vt_n-\vt_{n-1}),\zt^v_n\big)_{\FL} -F_{TR}(U_n,U_{n-1},\zt^v_n) + k A_p(U_n,\zt^v_n)\\
    &\qquad +k\theta A_F(U_n,\zt^v_n) + k(1-\theta) A_F(U_{n-1},\zt^v_n)
    -  (q_n,\zt^{v}_{n})_{\Gamma_q}\\
    &\qquad +k A_{div}(U_n,\zt^p_n) 
    +k A_{ALE}(U_n,\zt^{uf}_n)
    \Big\}\\
    &+\big(\ut(0)-\ut_0,\zt^{us}_0\big)_{\SO}
    +\big(\vt(0)-\vt_0,\zt^{v}_0\big)_{\FL}
    +\big(\vt(0)-\vt_0,\zt^{v}_0\big)_{\SO}
  \end{aligned}
\end{equation}

If the triplet $U_n=(p_n,\vt_n,\ut_n)\in X_h$ is the solution of the discrete fluid-structure interaction system of \eqref{aletime:2} in every time step $n=0,...,N$ with the control parameter $(q_n)_{n=1}^N$ in the boundary condition, the useful identity 
\begin{align}
j((q_n)_{n=1}^N):=J((q_n)_{n=1}^N,(U_n(q_n))_{n=0}^N)=L((q_n)_{n=1}^N,(U_n)_{n=0}^N,(Z_n)_{n=0}^N) 
\end{align}
is true for arbitrary values $(Z_n)\in Y_h$,  $n=0,\dots,N$. If we denote by $(\delta U_n)_{n=0}^N=\frac{d}{dq}(U_n)_{n=0}^N ((\delta q)_{n=1}^N)$ the derivative of the state variable with respect to the control, we obtain via the Lagrange functional the representation
\begin{multline*}
  j'((q_n)_{n=1}^N)((\delta q)_{n=1}^N)
  =L'_q((q_n)_{n=1}^N,(U_n)_{n=0}^N,(Z_n)_{n=0}^N)((\delta
  q)_{n=1}^N)\\
  +L'_{U}((q_n)_{n=1}^N,(U_n)_{n=0}^N,(Z_n)_{n=0}^N)(\delta
  U)_{n=0}^N) 
\end{multline*}
of the derivative of the reduced functional. If we choose the Lagrange multiplier $Z_n\in Y_h$ for $ n=N,...,0$ such that the dual problem
\begin{align}\label{eq_dual}
    \frac{d}{d U_n}  L((q_n)_{n=1}^N,(U_n)_{n=0}^N,(Z_n)_{n=0}^N)(\Phi)&=0 \quad \forall \Phi \in X_h \quad \text{for } n=N,\dots,0 ,
   \end{align}
   is fulfilled, then we can evaluate the derivative of the reduced functional $j((q_n)_{n=1}^N)$ in an arbitrary direction $(\delta q)_{n=1}^N$ by evaluating 
  \begin{align}\label{eq_gradient}
    j'((q_n)_{n=1}^N))((\delta q)_{n=1}^N)=L'_q((q_n)_{n=1}^N),(U_n)_{n=0}^N,(Z_n)_{n=0}^N)((\delta q)_{n=1}^N).
  \end{align}
 This enables us to apply any gradient based optimization algorithm. We will later use a limited memory version of the Broyden-Fletcher–Goldfarb-Shanno (BFGS) update formula (see for example \cite{Geiger2013}) to find a local minima of the discretized optimization problem.
 
 In every update step we first have to solve for the solution $U_n$ of the state equation \eqref{aletime:2} for $n=0,...,N$ and then compute the dual problem \eqref{eq_dual} for $n=N,...,0$. Thereby, the dual problem for $n=N-1,\dots,1$ consists of three equations with the test function $\Phi \in X_h$. Due to derivatives with respect to the velocity variable $v_n$, we obtain:
 \begin{multline}\label{dual:v}
   \big(\rho^0_s \phi,\zt^v_n\big)_{\SO}+k\theta \frac{d}{d v_n} A_S(U_n,\zt^v_n)(\phi) -k\theta\big(\phi,\zt^{us}_n\big)_\SO  
   +(\bar J_n \phi,\zt^v_n \big)_{\FL}
   \\-\frac{d}{d v_n} F_{TR}(U_n,U_{n-1},\zt^v_n)(\phi) 
     +k\theta \frac{d}{d v_n} A_F(U_n,\zt^v_n)(\phi) 
    +k \frac{d}{d v_n} A_{div}(U_n,\zt^p_n) (\phi)\\
    =k (\vt_n-\tilde \vt(t_n),\phi)_{\FL} + 
    \big(\rho^0_s \phi,\zt^v_{n+1} \big)_{\SO}
    \\ - k(1-\theta) \frac{d}{d v_n}A_S(U_{n},\zt^v_{n+1})(\phi)     -k(1-\theta)\big(\phi,\zt^{us}_{n+1}\big)_\SO    
    +(\bar J_{n+1} \phi,\zt^v_{n+1} \big)_{\FL}\\  +\frac{d}{d v_n} F_{TR}(U_{n+1},U_{n},\zt^v_{n+1})(\phi) 
    - k(1-\theta) \frac{d}{d v_n}A_F(U_{n},\zt^v_{n+1})(\phi).   
 \end{multline}
Due to derivatives of the Lagrangian with respect to the displacement $u_n$, we obtain:
\begin{multline}\label{dual:u}
  k\theta \frac{d}{d u_n} A_S(U_n,\zt^v_n)(\psi) +    \big(\psi,\zt^{us}_n\big)_{\SO}\\
  +( \frac{d}{d u_n}(\bar J_n)(\phi)  (\vt_n-\vt_{n-1}) ,\zt^v_n \big)_{\FL} 
  -\frac{d}{d u_n} F_{TR}(U_n,U_{n-1},\zt^v_n)(\psi) \\
  + k \frac{d}{d u_n} A_p(U_n,\zt^v_n)(\psi) +k\theta \frac{d}{d u_n} A_F(U_n,\zt^v_n)(\psi)\\
  +k \frac{d}{d u_n} A_{div}(U_n,\zt^p_n) (\psi)
  +k \frac{d}{d u_n} A_{ALE}(U_n,\zt^{uf}_n)(\psi)\\
  =
  - k(1-\theta) \frac{d}{d u_n}A_S(U_{n},\zt^v_{n+1})(\psi)  
  +\big(\psi,\zt^{us}_{n+1}\big)\\
  -( \frac{d}{d u_n}(\bar J_{n+1})(\phi)  (\vt_{n+1}-\vt_{n}) ,\zt^v_{n+1} \big)_{\FL} 
  +\frac{d}{d u_n} F_{TR}(U_{n+1},U_{n},\zt^v_{n+1})(\psi) \\
  - k(1-\theta) \frac{d}{d u_n}A_F(U_{n},\zt^v_{n+1})(\psi)
  +k (\ut_n-\tilde \ut(t_n),\psi)_{\SO}.
\end{multline}
Finally due to derivatives of the Lagrangian with respect to the pressure variable $p_n$, we obtain:
   \begin{equation}\label{dual:p}
     k \frac{d}{d p_n} A_p(U_n,\zt^v_n)(\xi)=0 .
   \end{equation}
The first and last step of the discrete dual problem have a slightly
different structure, but can be derived in a similar
way. {Since the monolithic formulation is a Petrov
  Galerkin formulation with different trial and test spaces, the
  adjoint coupling conditions differ from the primal ones. In the primal problem the solid displacement field enters as Dirichlet condition on the interface for the ALE extension problem. In the adjoint problem the shape derivatives of the adjoint ALE equation are coupled with the adjoint solid problem via a global test function which corresponds to a Neumann condition. As  $\zt^{uf}$ fulfills  zero Dirichlet conditions on the interface, this corresponds to a back coupling of the shape derivatives into the adjoint solid problem via residuum terms. Similar to the primal problem the adjoint velocity $\zt^v$ has to match on the interface and in addition an ``adjoint dynamic'' coupling condition is hidden in the test function $\phi$. Therefore, block preconditioners suggested in the literature cannot be directly applied to the adjoint problem, but have to be adapted to the new structure.}

\subsection{Short notation for state and dual equation}
\label{sec:red}

\paragraph{Short notation of the state equation} 
Key to the efficiency of the multigrid approach demonstrated
in~\cite{FailerRichter2019} is a condensation of the deformation
unknown $\ut_n$ from the solid problem. 
The last equation in~(\ref{aletime:2}) gives a relation for the new
deformation at time~$t_n$
\begin{align}
\ut_n = \ut_{n-1}+k\theta \vt_n + k(1-\theta)\vt_{n-1}\text{ in }\SO 
\label{eq_vel_disp_rel}
\end{align}
and we will use this representation to eliminate the unknown
deformation and base the 
solid stresses purely on the last time step and the unknown velocity, i.e.
by expressing the deformation gradient as 
\begin{equation}\label{dispvel}
  \begin{aligned}
    \Ft_n=\Ft(\ut_{n}) \,&\widehat{=}\, \Ft(\ut_{n-1},\vt_{n-1};\vt_n) \\
    &=I+\nabla   \big(\ut_{n-1}+k\theta \vt_n + 
    k(1-\theta)\vt_{n-1}\big)\text{ in }\SO. 
  \end{aligned}
\end{equation}
Removing the solid deformation from the momentum equation will help to
reduce the algebraic systems in Section~\ref{sec:solver}. A similar
technique within an Eulerian formulation and using a characteristics
method is presented in~\cite{Pironneau2016,Pironneau2019}.

   For each time step $t_{n-1}\mapsto t_n$ we introduce the following
short notation for the system of algebraic equations that is based on
the splitting of the solution into unknowns acting in the fluid domain
$(\vt_f,\ut_f)$, on the interface $(\vt_i,\ut_i)$  and those on the solid
$(\vt_{s},\ut_{s})$.  The pressure variable $p$ acts in the fluid
and on the interface.
\begin{equation}\label{system}
  \underbrace{\begin{pmatrix}
    {\cal D}(p,\vt_f,\ut_f,\vt_{i},\ut_{i},\vt_{s},\ut_{s}) \\    
    {\cal M}^f(p,\vt_f,\ut_f,\vt_{i},\ut_i) \\
    {\cal M}^i(p,\vt_f,\ut_f,\vt_{i},\ut_i,\vt_s) \\
    {\cal M}^s(p,\vt_{i},\ut_i,\vt_s) \\
    {\cal E}(\ut_f,\ut_i)\\
    {\cal U}^i(\vt_{i},\ut_{i},\vt_s,\ut_s)\\
    {\cal U}^s(\vt_{i},\ut_{i},\vt_s,\ut_s)
  \end{pmatrix}}_{=:  {\cal A}(U)}
  =
  \underbrace{
    \begin{pmatrix}
      {\cal B}_1 \\ {\cal B}_2\\ {\cal B}_3\\ {\cal B}_4\\
      {\cal B}_5\\ {\cal B}_6 \\{\cal B}_7
    \end{pmatrix}}_{=: {\cal B}}
\end{equation}
${\cal D}$ describes the divergence equation which acts in the
fluid domain and on the interface, ${\cal M}$
the two momentum equations, acting in the fluid domain, on the
interface and in the solid domain (which is indicated by a
corresponding index), ${\cal E}$ describes the ALE extension in the fluid domain
and ${\cal U}$ is the relation between solid velocity and
solid deformation, acting on the interface degrees of freedom and in
the solid. Note that ${\cal M}^i$ and ${\cal M}^s$, the term describing the
momentum equations, do not directly depend on the solid deformation 
$\ut_{s}$ as we base the deformation gradient on the velocity,
see~(\ref{dispvel}).

\paragraph{Short notation Dual Equation}
We aim at applying a similar reduction scheme to the adjoint
problem. Here, there is no direct counterpart
to~(\ref{eq_vel_disp_rel}). Instead, we first introduce the new
variable $\tilde \zt^{us}_n$ such that 
\begin{align}\label{eq_sub_zus}
\big(\psi,\tilde \zt^{us}_n\big)_{\SO}=  \big(\psi, \zt^{us}_n\big)_{\SO}+k\theta \frac{d}{d u_n} A_S(U_n,\zt^v_n)(\psi).
\end{align}
{Thereby we can substitute all terms in~(\ref{dual:v}),~(\ref{dual:u}) and~(\ref{dual:p}) which depend on $\zt^{us}_n$ by the new variable $\tilde \zt^{us}_n$, such as}
\begin{align}
-\theta k \big(\phi,\zt^{us}_n\big)_{\SO}=- \theta k \big(\phi,\tilde \zt^{us}_n\big)_{\SO}+(\theta k)^2\frac{d}{d u_n} A_S(U_n,\zt^v_n)(\phi)
\end{align}
in~(\ref{dual:v}). Now the adjoint terms ${\cal M}_{\ut_i}$ and ${\cal M}_{\ut_s}$ resulting from derivatives of the momentum equation with respect to the displacement variable do not depend on the adjoint velocity variable $\zt^{v}$ anymore which will enable later to decouple the problem in three well conditioned subproblems. Furthermore the "adjoint dynamic`` coupling conditions now corresponds to equivalents of adjoint boundary forces on the interface as in the state equation.
 
For each time step $t_{n+1}\mapsto t_n$ we introduce again a
short notation for the system of algebraic equations that is based on
the splitting of the adjoint solution into unknowns acting in the fluid domain
$(\zt^v_f,\zt^{uf}_f)$, on the interface $(\zt^v_i, \tilde \zt^{us}_i)$  and those on the solid
$(\zt^v_s, \tilde \zt^{us}_s)$.  The adjoint pressure variable $z^p$ acts in the fluid
and on the interface. 
   \begin{equation}\label{system}
     \underbrace{\begin{pmatrix}
     {\cal M}_{p} (\zt^v_f,\zt^v_{i},\zt^v_{s}) \\
    {\cal D}_{\vt_f} (\zt^p)+ {\cal M}_{\vt_f} (\zt^v_f,\zt^v_i)\\
    {\cal D}_{\ut_f} (\zt^p)+ {\cal M}_{\ut_f} (\zt^v_f,\zt^v_i)+{\cal E}_{\ut_f}(\zt^{uf}_f)\\    
    {\cal D}_{\vt_i} (\zt^p)+ {\cal M}_{\vt_i} (\zt^v_f,\zt^v_i,\zt^v_s) +{\cal U}_{\vt_i} ( \tilde \zt^{us}_{i}, \tilde\zt^{us}_{s}) \\
    {\cal D}_{\ut_i} (\zt^p)+ {\cal M}_{\ut_i} (\zt^v_f,\zt^v_i)+{\cal E}_{\ut_i}(\zt^{uf}_f)+{\cal U}_{\ut_i} (\tilde \zt^{us}_{i}, \tilde\zt^{us}_{s})\\
    {\cal M}_{\vt_s} (\zt^v_i,\zt^v_s)+{\cal U}_{\vt_s} ( \tilde \zt^{us}_{i}, \tilde \zt^{us}_{s}) \\
    {\cal U}_{\ut_s} ( \tilde \zt^{us}_{i}, \tilde \zt^{us}_{s})\\
     \end{pmatrix}}_{=:
       {{\cal A}^{\text{Dual}}}(Z)}
     =
     \underbrace{\begin{pmatrix}
         {\cal B}_1^d \\ {\cal B}_2^d\\ {\cal B}_3^d\\ {\cal B}_4^d\\
         {\cal B}_5^d\\ {\cal B}_6^d \\{\cal B}_7^d
       \end{pmatrix}}_{=: {\cal B}^d}
\end{equation}

${\cal M}_{p} $ describes the adjoint divergence equation which acts in the
fluid domain and on the interface, ${\cal M}_{\vt}$ and ${\cal M}_{\ut}$
the derivatives of the momentum equation with respect to the velocity and displacement variable, acting in the fluid domain, on the
interface and in the solid domain (which is indicated by a
corresponding index) and ${\cal E}_{\ut}$ describes the adjoint ALE extension in the fluid domain
and ${\cal U}_{\vt} $ and ${\cal U}_{\ut} $ result from the relation between solid velocity and
solid deformation, which act on the interface degrees of freedom and in
the solid.

\section{Solution of the algebraic systems}\label{sec:solver}

In~\cite{FailerRichter2019}, we have derived an efficient approximated
Newton scheme for the forward fluid-structure interaction 
problem. We briefly outline the main steps and then focus on
transferring these ideas to the  dual equations. 
The general idea is described by the following two steps
\begin{enumerate}
\item In the Jacobian, we omit the derivatives of the Navier-Stokes
  equations with respect to the fluid domain
  deformation, which results in an approximated Newton
  scheme. In~\cite[chapter 5]{Richter2017} it is documented that this
  approximation will slightly increase the iteration counts of the
  Newton scheme. On the other hand, the overall computational time is
  nevertheless reduced, since assembly times for these neglected terms are
  especially high. Since the Newton residual is not changed, the
  resulting nonlinear solver is of an approximated Newton type. 
\item We use the discretization of the relation $\partial_t \ut=\vt$
  between solid deformation  and solid velocity, namely
  $\ut^{n+1}=\ut^n + \theta k \vt^{n+1}+(1-\theta)k\vt^n$ to
  reformulate the solid's deformation gradient based on the velocity
  instead of the deformation. This step has been explained in the
  previous section. 
\end{enumerate}
These two steps, the first one being an approximation, while the
second is an equivalence transformation, allow to reduce the number of
couplings in the Jacobian in such a way that each linear step
falls apart into three successive linear systems. The first one
describes the coupled momentum equation for fluid- and
solid-velocity, the second realizes the solid's velocity-deformation
relation and the third one stands for the ALE extension.
We finally
note that the approximations only involve the Jacobian. The residual
of the systems is not altered such that we still solve the original
problem and compute the exact discrete gradient.

{Then, in Section~\ref{solve:dual} we describe the extension of this solution mechanism to the adjoint system. Two major differences occur: first, the adjoint system is linear, such that we realize the solver in the framework of a preconditioned Richardson iteration. The preconditioner takes the place of the approximated Jacobian. Second, the adjoint interface coupling conditions differ from the primal conditions as outlined in the last paragraph of Section~\ref{Sec::OptSystem}. This will call for a modification of the condensation procedure introduced as second reduction step in the primal solver. 
}
\subsection{Solution of the primal problem}

In each time step of the forward problem we must solve a nonlinear
problem. We employ an approximated Newton scheme
\begin{equation}\label{newton}
  \tilde {\cal A}'(U^{(l)}) W^{(l)} = {\cal B}-{\cal A}(U^{(l)}),\quad
  U^{(l+1)}:=U^{(l)}+\omega^{(l)} \cdot W^{(l)},
\end{equation}
where $\omega^{(l)}$ is a line search parameter, $U^{(0)}$ an initial
guess. By ${\cal A}'(U)$ we denote the Jacobian, by $\tilde {\cal
  A}'(U)$ an approximation. As outlined in~\cite{FailerRichter2019}
the Jacobian is modified in two essential steps: first, in the
Navier-Stokes problem, we skip the derivatives with respect to the ALE
discretization. These terms are computationally expensive and they
further introduce the only couplings from the fluid problem to the
deformation unknowns. In~\cite[chapter 5]{Richter2017} it has been
shown that while this approximation does slightly worsen Newton's
convergence rate, the overall efficiency is nevertheless increased, as
the number of additional Newton steps is very small in comparison to
the savings in assembly time. Second, we employ the reduction step
outlines in Section~\ref{sec:red}, which is a static
condensation of the deformation unknowns from the solid's momentum
equation. Taken together, both steps completely remove all deformation
couplings from the combined fluid-solid momentum equation and the
Jacobian takes the form
\begin{equation}\label{REDUCEDMATRIX}
  \left(
  \begin{array}{cccc|c|cc}
    0            & {\cal D}_{\vt_f} & {\cal D}_{\vt_i} & 0 &  
    0 &  0&0\\ 
    {\cal M}^f_p & {\cal M}^f_{\vt_f} & {\cal M}^f_{\vt_i} & 0 &  0 &  0&0 \\
    {\cal M}^i_p & {\cal M}^i_{\vt_f} & {\cal M}^i_{\vt_i} & {\cal M}^i_{\vt_s} &  0 &  0&0\\
    {\cal M}^s_p & 0  & {\cal M}^s_{\vt_i} & {\cal M}^s_{\vt_s} & 0 & 0&0 \\
    \hline
    0 & 0 & 0 & 0 & {\cal E}^f_{\ut_{f}} & {\cal E}^f_{\ut_{i}}&0\\
    \hline
    0&0&{\cal U}^i_{\vt_{i}}&{\cal U}^i_{\vt_{s}}&0&{\cal U}^i_{\ut_{i}}&{\cal U}^i_{\ut_{i}}\\
    0&0&{\cal U}^s_{\vt_{i}}&{\cal U}^s_{\vt_{s}}&0&{\cal U}^s_{\ut_{i}}&{\cal U}^s_{\ut_{i}}\\
  \end{array}
  \right)
  \left(
  \begin{array}{c}
    \delta \pt  \\
    \deltat \vt_f \\
    \deltat \vt_{i} \\
    \deltat \vt_{s} \\
    \deltat \ut_f \\ 
    \deltat \ut_{i}\\
    \deltat \ut_{s}
  \end{array}
  \right) = \left(
  \begin{array}{c}
    \bt_1 \\ \bt_2 \\ \bt_3 \\ \bt_4 \\ \bt_5 \\ \bt_6  \\ \bt_7
  \end{array}
  \right).
\end{equation}
This reduced linear system decomposes into three sub-steps. First, the
coupled momentum equation, living in fluid and solid domain and acting
on pressure and velocity,
   \begin{equation}\label{problem:1}
     \left(
     \begin{array}{cccc}
       0            & {\cal D}_{\vt_f} & {\cal D}_{\vt_i} & 0\\
       {\cal M}^f_p & {\cal M}^f_{\vt_f} & {\cal M}^f_{\vt_i} & 0\\
       {\cal M}^i_p & {\cal M}^i_{\vt_f} & {\cal M}^i_{\vt_i} & {\cal M}^i_{\vt_s}\\
       {\cal M}^s_p & 0  & {\cal M}^s_{\vt_i} & {\cal M}^s_{\vt_s}\\
     \end{array}
     \right)
     \left(
     \begin{array}{c}
       \delta \pt  \\
       \deltat \vt_f \\
       \deltat \vt_{i} \\
       \deltat \vt_{s} \\
     \end{array}
     \right) = \left(
     \begin{array}{c}
       \bt_1 \\ \bt_2 \\ \bt_3 \\ \bt_4
     \end{array}
     \right),
   \end{equation}
second, the update equation for the deformation on the interface and
within the solid domain,
    \begin{equation}\label{problem:2}
      \left(
      \begin{array}{cc}
        {\cal U}^i_{\ut_{i}}&{\cal U}^i_{\ut_{i}}\\
        {\cal U}^s_{\ut_{i}}&{\cal U}^s_{\ut_{i}}\\
      \end{array}
      \right)
      \left(
      \begin{array}{c}
        \deltat \ut_{i}\\
        \deltat \ut_{s}
      \end{array}
      \right) = \left(
      \begin{array}{c}
        \bt_6  \\ \bt_7
      \end{array}
      \right)-
      \begin{pmatrix}
        {\cal U}^i_{\vt_{i}}&{\cal U}^i_{\vt_{s}}\\
        {\cal U}^s_{\vt_{i}}&{\cal U}^s_{\vt_{s}}
      \end{pmatrix}
      \begin{pmatrix} \vt_i \\ \vt_s \end{pmatrix},
    \end{equation}
which is a finite element discretization of the zero-order equation 
$\ut_n = \ut_{n+1}+k (1-\theta)\vt_{n-1}+k\theta\vt_n$. This update
can be performed by one algebraic vector-addition. Finally it remains
to solve for the ALE extension equation
\begin{equation}\label{problem:3}
  {\cal E}^f_{\ut_{f}}\delta\ut_{f} = \bt_5- {\cal E}^f_{\ut_{i}}\delta\ut_{i},
\end{equation}
one simple equation, usually either a vector Laplacian or a linear
elasticity problem, see~\cite[Section 5.2.5]{Richter2017}.

\subsection{Dual}\label{solve:dual}

Due to the unsymmetrical structure of the fluid-structure interaction
model the block collocation and coupling of the blocks in the
transposed Jacobian $A'(U)^T$ in the dual problem differs to the
Jacobian of the primal problem. This stays in strong relation to the
adjoint coupling conditions, see Section \ref{Sec::OptSystem}. Hence,
block preconditioners developed for the state problem can not be
applied in a black box way to the linear systems arising in the dual
problem, but have to be adjusted. Furthermore, the dual system is
linear such that the approximated Newton scheme must be replaced by a
different concept. We start by indicating the full system matrix of the
dual problem
\begin{equation}\label{DUALFULLMATRIX}
  \underbrace{\left(
    \begin{array}{cccc|c|cc}
      0&{\cal M}^{f,T}_p &  {\cal M}^{i,T}_p & {\cal M}^{s,T}_p &0&0&0\\
      {\cal D}_{\vt_f}^T &{\cal M}^{f,T}_{\vt_f} &{\cal M}^{i,T}_{\vt_f} &0&0&0&0\\
      \SCC {\cal D}_{\ut_f}^T &\SCC {\cal M}^{f,T}_{\ut_f} &\SCC {\cal
        M}^{i,T}_{\ut_f} &0&{\cal E}^{f,T}_{\ut_{f}} &0&0\\
      \hline
          {\cal D}^T_{\vt_i} &{\cal M}^{f,T}_{\vt_i} &{\cal M}^{i,T}_{\vt_i} &{\cal M}^{s,T}_{\vt_i} &0&{\cal U}^{i,T}_{\vt_{i}}&{\cal U}^{s,T}_{\vt_{i}}\\
          \SCC {\cal D}_{\ut_i}^T&\SCC {\cal M}^{f,T}_{\ut_i}&\boldsymbol{ \SCC
            {\cal M}^{i,T}_{\ut_i}}&\boldsymbol{ {\cal     M}^{s,T}_{\ut_i}}& {\cal
            E}^{f,T}_{\ut_{i}}&{\cal U}^{i,T}_{\ut_{i}}&{\cal U}^{s,T}_{\ut_{i}}\\
          \hline
          0&0&{\cal M}^{i,T}_{\vt_s}  &{\cal M}^{s,T}_{\vt_s}&0&{\cal U}^{i,T}_{\vt_{s}}&{\cal U}^{s,T}_{\vt_{s}}\\
          0&0&\boldsymbol{{\cal M}^{i,T}_{\ut_s}}&\boldsymbol{{\cal M}^{s,T}_{\ut_s}}&0&{\cal U}^{i,T}_{\ut_{s}}&{\cal U}^{s,T}_{\ut_{s}}\\
    \end{array}
    \right)}_{=A'_D}
    \left(
    \begin{array}{c}
      \zt^p  \\
      \zt^v_f \\ \zt^v_i \\ 
      \zt^v_s \\ \zt^{uf}_f\\
      \tilde \zt^{us}_i \\ \tilde  \zt^{us}_s
    \end{array}
    \right) = \left(
    \begin{array}{c}
     {\cal B}^d_1 \\ {\cal B}^d_2 \\ {\cal B}^d_3 \\ {\cal B}^d_4 \\ {\cal B}^d_5 \\ {\cal B}^d_6  \\ {\cal B}^d_7
    \end{array}
    \right),
\end{equation}
given as the transposed of the primal Jacobian, $A_D =
A'(U)^T$, see~\cite{FailerRichter2019}. 

For solving the dual problem we want to mimic the primal approach:
first, approximate the system matrix by neglecting couplings, second,
use the static condensation as described in Section~\ref{sec:red}
in~\eqref{eq_sub_zus}. As the problem is linear, a direct modification
of the system matrix would alter the dual solution. Instead, we
approximate the solution by a preconditioned Richardson
iteration with an inexact matrix $\tilde A'_D\approx A'(U)^T$  as
preconditioner (approximated by a geometric multigrid solver)
\begin{equation}\label{dualiteration}
  Z^{(0)}=0,\quad
  \tilde A'_D W^{(l)}={\cal B}^d-A'(U)^TZ^{(l-1)},\quad
  Z^{(l)}=Z^{(l-1)}+W^{(l)},
\end{equation}
where $Z^{(l)} =
\{\zt^p,\zt^v_f,\zt^v_i,\zt^v_s,\zt^{uf}_f,\tilde \zt^{us}_i,\tilde \zt^{us}_s\}$
and the update in every Richardson iteration is given by $W^{(l)} =
\{\delta\zt^p,\delta\zt^v_f,\delta\zt^v_i,\delta\zt^v_s,\delta\zt^{uf}_f,\delta \tilde \zt^{us}_i,\delta \tilde \zt^{us}_s\}$. 
The residual
is computed based on the full Jacobian $A'(U)^T$ (including the
ALE derivatives) such that we still converge to the original adjoint
problem. Since we never assemble the complete Jacobian $A'(U)$ in the
primal solver, the adjoint residual ${\cal B}^d-A'(U)^TZ^{(l-1)}$ can
be established in a matrix free setting. 

Then, similar to the described approach in case of the primal system,
we neglect the ALE terms (shaded entries). Finally, we reorder to
reach the preconditioned iteration
\begin{equation}\label{DUALREDUCEDMATRIX}
  \underbrace{\left(
  \begin{array}{cccc|c|cc}
    0&{\cal M}^{f,T}_p &  {\cal M}^{i,T}_p & {\cal M}^{s,T}_p &0&0&0\\
    {\cal D}_{\vt_f}^T &{\cal M}^{f,T}_{\vt_f} &{\cal M}^{i,T}_{\vt_f} &0&0&0&0\\
    {\cal D}^T_{\vt_i} &{\cal M}^{f,T}_{\vt_i} &{\cal M}^{i,T}_{\vt_i} &{\cal M}^{s,T}_{\vt_i} &0&{\cal U}^{i,T}_{\vt_{i}}&{\cal U}^{s,T}_{\vt_{i}}\\
    0&0&{\cal M}^{i,T}_{\vt_s}  &{\cal M}^{s,T}_{\vt_s}&0&{\cal U}^{i,T}_{\vt_{s}}&{\cal U}^{s,T}_{\vt_{s}}\\
    \hline
    0 &0 &0 &0&{\cal E}^{f,T}_{\ut_{f}} &0&0\\
    \hline
    0&0&0&0& {\cal   E}^{f,T}_{\ut_{i}}&{\cal U}^{i,T}_{\ut_{i}}&{\cal U}^{s,T}_{\ut_{i}}\\
    0&0&0&0&0&{\cal U}^{i,T}_{\ut_{s}}&{\cal U}^{s,T}_{\ut_{s}}\\
  \end{array}
  \right)}_{=\tilde A'_D}
      \left(
    \begin{array}{c}
      \delta \zt^p  \\
      \deltat \zt^v_f \\ \deltat \zt^v_i \\ 
      \deltat \zt^v_s \\ \deltat \zt^{uf}_f\\
      \deltat \tilde \zt^{us}_i \\ \deltat \tilde \zt^{us}_s
    \end{array}
    \right) 
  = 
  {\cal B}^d-A'(U)^T  Z^{(l-1)},
  \end{equation}
with the preconditioner $\tilde A_D$ that decomposes into three
separate steps. First, the equation for the adjoint mesh deformation
variable  
\begin{equation}\label{equ:dualproblem1}
  {\cal E}^{f,T}_{\ut_f} \deltat \zt^{uf}_f=\bt^d_3. 
\end{equation}
Usually a symmetric extension operator ${\cal E}^f_{\ut_f}$ can be
chosen. This avoids re-assembly of this matrix and possible
preparations for the linear solver. See~\cite[Section
  5.3.5]{Richter2017} for different efficient options for extension
operators. Second, the update for the adjoint solid deformation, 
\begin{equation}\label{equ:dualproblem2}
  \left(
  \begin{array}{cc}
    {\cal U}^{i,T}_{\ut_{i}}&{\cal U}^{s,T}_{\ut_{i}}\\
    {\cal U}^{i,T}_{\ut_{s}}&{\cal U}^{s,T}_{\ut_{s}}\\
  \end{array}
  \right)
  \left(
  \begin{array}{c}
     \deltat \tilde \zt^{us}_i\\
     \deltat \tilde \zt^{us}_s\\
  \end{array}
  \right) = \left(
  \begin{array}{c}
   \bt^d_5 -{\cal E}^{f,T}_{\ut_i} \deltat \zt^{uf}_f \\
   \bt^d_7
  \end{array}
  \right),
\end{equation}
which only involves inversion of the mass matrix and finally the
update for the adjoint velocity and adjoint pressure 
\begin{equation}\label{equ:dualproblem3}
  \left(
  \begin{array}{cccc}
    0&{\cal M}^{f,T}_p &  {\cal M}^{i,T}_p & {\cal M}^{s,T}_p\\
    {\cal D}_{\vt_f}^T &{\cal M}^{f,T}_{\vt_f} &{\cal M}^{i,T}_{\vt_f}&0\\
    {\cal D}^T_{\vt_i} &{\cal M}^{f,T}_{\vt_i} &{\cal M}^{i,T}_{\vt_i} &{\cal M}^{s,T}_{\vt_i}\\
    0&0&{\cal M}^{i,T}_{\vt_s}  &{\cal M}^{s,T}_{\vt_s}
  \end{array}
  \right)
  \left(
  \begin{array}{c}
    \delta \zt^p  \\
    \deltat \zt^v_f\\ 
    \deltat \zt^v_i \\ 
    \deltat \zt^v_s\\ 
  \end{array}
  \right) = \left(
  \begin{array}{c}
    \bt^d_1 \\ \bt^d_2 \\ \bt^d_4 \\ \bt^d_6
  \end{array}
  \right)-
    \left(
  \begin{array}{cc} 
    0&0\\
    0&0\\
    {\cal U}^{i,T}_{\vt_s}&{\cal U}^{s,T}_{\vt_s}\\
    {\cal U}^{i,T}_{\vt_s}&{\cal U}^{s,T}_{\vt_s}\\
  \end{array}
  \right)
  \left(
  \begin{array}{c}
     \deltat \tilde \zt^{us}_i\\
     \deltat \tilde \zt^{us}_s
  \end{array}
  \right).
\end{equation}

The numbering of the right hand side $\bt^d_1,\dots,\bt^d_7$ is according
to~(\ref{DUALFULLMATRIX}). 
As we do not modify the residuum, the derivatives with respect to the ALE transformation ${\cal M}_{\ut_f} (\zt^v_f,\zt^v_i)$ and  ${\cal M}_{\ut_i} (\zt^v_f,\zt^v_i)$ still enter into $\bt_3$ and $\bt_5$. 
Hence the resulting problem in Equation~\eqref{equ:dualproblem1}
corresponds to a linear elasticity problem on the fluid domain with an
artificial forcing term in the right hand side and zero Dirichlet data
on the interface. In Equation~\eqref{equ:dualproblem2} the shape
derivatives of the the ALE transformation enter via Residuum terms
${\cal M}_{\ut_i} (\zt^v_f,\zt^v_i)$ on the interface. These terms
contain the adjoint geometric coupling condition. An explicit update
by one vector-addition as for the corresponding primal equation is not
possible. The ``adjoint kinematic'' and ``adjoint dynamic'' coupling
conditions are fully incorporated in \eqref{equ:dualproblem3}, similar
as for the state equation, and thereby these
coupling conditions are fully resolved in every Richardson iteration.

\subsection{Solution of the linear problems}

In each step of the Newton iteration for solving the state equation
and in each step of the Richardson iteration in the case of the adjoint system, we
must approximate three individual linear systems of equations. The
mesh-update problems
are usually of elliptic type, the vector Laplacian or a linear
elasticity problem. Here, standard geometric multigrid solvers are
highly efficient. Problem~(\ref{equ:dualproblem2}) and the primal
counterpart correspond to zero order
equations. Multigrid solvers or the CG method converge with optimal
efficiency. It remains to approximate the coupled momentum
equations, given by~(\ref{equ:dualproblem3}) in the dual case. Here,
we are lacking any desirable structure. The matrices are not
symmetric, they 
feature a saddle-point structure and involve different scaling of the
fluid- and solid-problem. We approximate these equations by a GMRES
iteration that is preconditioned with a geometric multigrid
solver. Within the multigrid iteration we employ a smoother of Vanka
type, where we invert local patches exactly. These patches correspond
to all degrees of freedom of one element (within the fluid) and to a
union of $2^d$ elements (within the solid). For the sake of
simplicity (in terms of implementational effort) we use this highly
robust solver also for the other two problems, despite their simpler
character. For details we refer to~\cite{FailerRichter2019}.

{\subsection{Algorithm}
To get an overview how the final optimization routine works we summarized all the intermediate steps in the following algorithm: 
\begin{algorithm}
  \caption{Optimization loop}
  {
  \begin{algorithmic}
	\\	
	\State Set $(q_n)_{n=1}^N\in \R^N$, $U_0=U(0)$ 
	\While{ $\norm{\nabla j((q_n)_{n=1}^N)}> tol_{opt}$}
		\State Set $q=q_{old}+ \alpha \cdot d_n$ 
		\State {Compute the solution $(U_n)_{n=0}^N$ of the primal problem \eqref{aletime:2}}:
		\For {$n=1:N$}
		\While{$\rho_n>tol_{Newton}$} 
		\State Solve linear systems  \eqref{problem:1}, \eqref{problem:2}, \eqref{problem:3} (Preconditioned GMRES)
\State Apply Newton update (see \eqref{newton}) 
\State Compute Newton residuum $\rho_n$ 
		\EndWhile
		\EndFor
		\If{$j((q_n)_{n=1}^N)>j((q_{old,n})_{n=1}^N)$ }
	$\alpha=0.5\cdot\alpha$ and \textbf{continue}
\EndIf
		\State {Compute the solution $(Z_n)_{n=0}^N$  of the dual problem \eqref{eq_dual}:}
		\For {$n=N:1$}
		\While{$\rho_n>tol_{Richardson}$} 
			\State Solve linear system \eqref{equ:dualproblem1}, \eqref{equ:dualproblem2}, \eqref{equ:dualproblem3} (Preconditioned GMRES)
			\State Apply Richardson update (see \eqref{dualiteration}) 
		    \State Compute Richardson residuum $\rho_n$
		\EndWhile
		\EndFor
		\State {Evaluate gradient $\nabla j((q_n)_{n=1}^N)(\delta q)$ (see \eqref{eq_gradient})}
		\State Update BFGS Matrix and compute BFGS update $(d_n)_{n=1}^N$ 
	\EndWhile
\end{algorithmic}
}
\end{algorithm}
}

\section{Numerical Results}\label{sec:num}

\subsection{Problem configuration 2d}

  \begin{figure}[h]
  \centering
    \begin{tikzpicture}[scale=7 ]
	\draw (0,0) -- (0.2,0)  ;
	\draw (0.35,0) -- (1.5,0)  ;
	\draw (0.2,0) -- (0.2,-0.1)  ;
	\draw (0.2,-0.1) -- (0.35,-0.1)  node[black,anchor= north east ] {$\Gamma_{q}$};  ;
	\draw (0.35,-0.1) -- (0.35,0)  ;
	\draw (1.5,0) -- (1.5,0.41)  ;
	\draw (1.5,0.41) -- (0,0.41)  ;
	\draw (0,0) -- (0,0.41)  ;
	\draw (0.2,0.2) circle [radius=0.05];
	\draw (0.25,0.19) rectangle (0.6,0.21)  ;
	\draw (0,0.2) node[black,anchor= east ] {$\Gamma_f^{in}$};
	\draw (0.14,0.14) node[black ] {$\Gamma_f$};
	\draw (0.8,0.41) node[black ,anchor=south] {$\Gamma_f$};
	\draw (0.8,0.0) node[black ,anchor=north] {$\Gamma_f$};
	\draw (1.0,0.2) node[black] {$\FL$};
	\draw (1.5,0.2) node[black,anchor=west] {$\Gamma_f^{out}$};
	\draw (0.5,0.2)   node[black,anchor=north] {$\IN$};
	\draw [thick,->] (0,0.2) -- (0.1,0.2); 
	\draw [thick,->] (0,0.275) -- (0.08,0.275); 
	\draw [thick,->] (0,0.35) -- (0.04,0.35); 
	\draw [thick,->] (0,0.125) -- (0.08,0.125);
	\draw [thick,->] (0,0.05) -- (0.04,0.05); 
	\draw [thick,<-] (0.25,0.2) -- (0.35,0.3) node[black,anchor= west ] {$\Gamma_s$}; 
	\draw [thick,<-] (0.6,0.2) -- (0.7,0.3) node[black,anchor=south] {$A=(0.6,0.2)$} ; 
	\draw [thick,<-] (0.5625,0.2) -- (0.7,0.1) node[black,anchor=north] {$B=(0.5625,0.2)$} ; 
	\draw [thick,<-] (0.4,0.2) -- (0.3,0.1) node[black,anchor=north] {$\SO$}; 
	\draw [blue](0.525,0.19) rectangle (0.6,0.21) ;
\end{tikzpicture}
  \caption{Geometry for flow around cylinder with elastic beam. {The blue region denotes the observation domain $\Omega_{obs}$.}}\label{configuration_opt_fsi_2D}
  \end{figure}
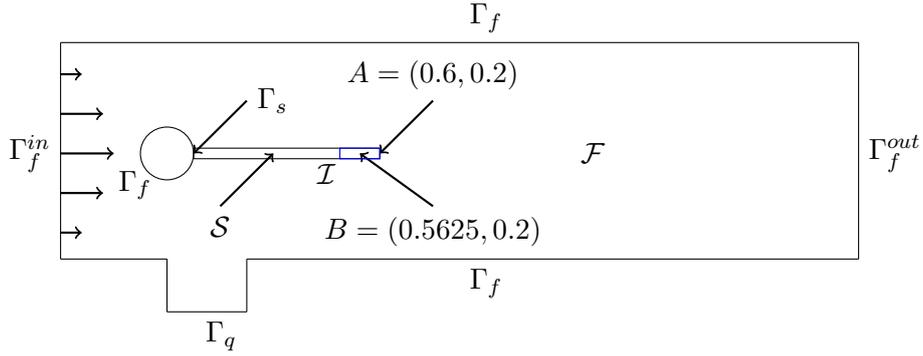
  
We modify the well known FSI Benchmark from Turek and Hron \cite{HronTurek2006} by adding an additional boundary $\Gamma_{q}$ as in Figure \ref{configuration_opt_fsi_2D}. The material parameters are chosen as for the FSI 1 Benchmark. 
In the solid the Lame parameters with $\lambda=2.0 \cdot 10^6\SI{}{\kg \per \m \s^2}$ and $\mu=0.5\cdot10^6\SI{}{ \kg \per \m \s^2}$ and a fluid viscosity $\nu_f=\SI{0.001}{\m^2 \per \s}$ are chosen. The solid and fluid densities are given by $\rho_s=\SI{1000.0}{kg\per \m^3}$ and $\rho_f=\SI{1000.0}{\kg \per \m^3}$. The inflow velocity is increased slowly during the time interval $I=[0,2]$ as for the instationary FSI 2 and FSI 3 benchmark.

In the first example, we would like to determine the pressure profile $q(t)\in L^2(I)$ at the control boundary  $\Gamma_{q}$ on the time interval $I=[0,6]$, leading to the displacement profile over time 
\begin{align}
\tilde u(t)=     
\begin{cases}
0 &\quad  t< \SI{2}{\s}\\
 0.01\cdot \sin(2\pi t) &\quad t \geq \SI{2}{\s}          
\end{cases}
\end{align}
in y-direction in the area $\Omega_{obs}=\{0.525\leq x\leq0.6,0.19\leq y\leq0.21 \}$ at the tip of the flag (see Fig.\ref{configuration_opt_fsi_2D}). To do so, we minimize the functional
\begin{equation}
\begin{aligned}
\min_{q \in  L^2(I)} J(q,\ut)= \frac{1}{2} \int_{0}^{6}\norm{\ut_y-\tilde u}^2_{\Omega_{Obs}}\text{ d}t+ \frac{\alpha}{2}\int_{0}^{6} q(t)^2 \text{ d}t
\end{aligned}
\end{equation}
constrained by the fluid-structure interaction problem. We discretize the system as presented in Section \ref{sec:disc} in time using a shifted Crank-Nicolson time stepping scheme with time step size $k=0.01$ and $\theta=0.5+2k$. The control variable is chosen to be piece-wise constant on every time interval ($\dim(Q)= 600$). The Tikhonov regularization parameter is set to $\alpha=1.0\cdot10^{-17}$.

\subsection{Problem configuration 3d}
\begin{figure}
  \centering
  \includegraphics[scale=0.2]{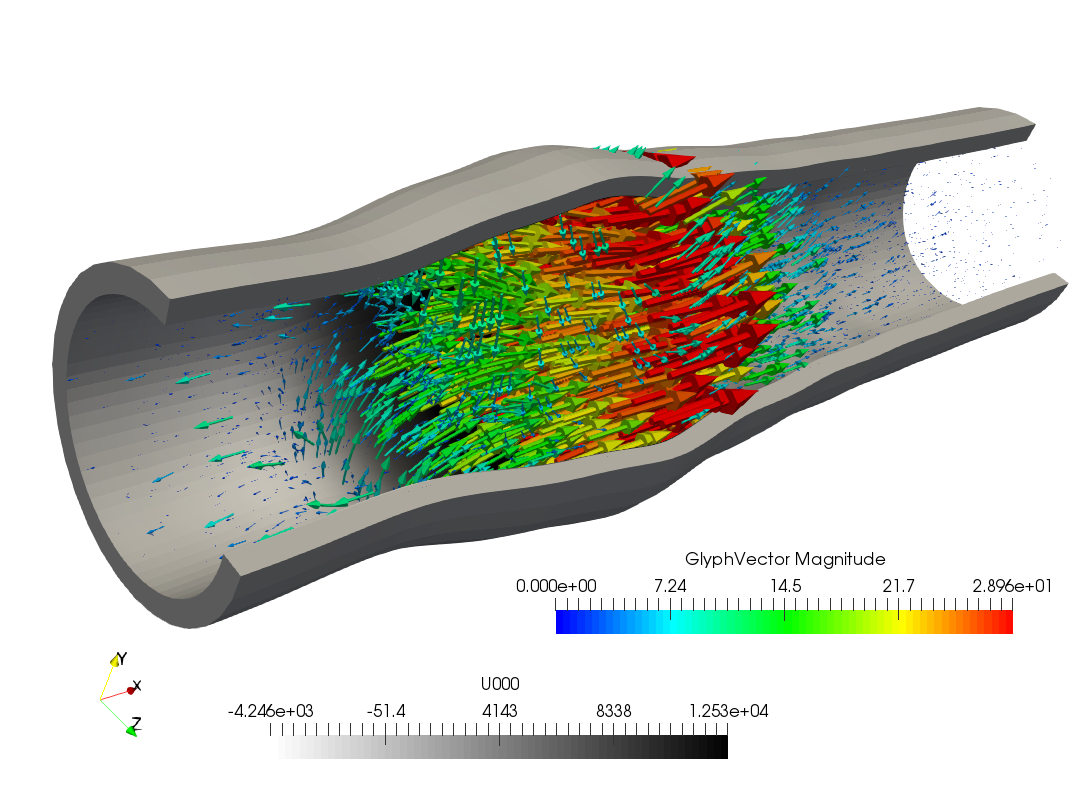}
  \caption{Velocity field of the pressure wave at $t=\SI{0.006}{\s}$ on the deformed domain (amplified by a factor 10) }\label{Solution_3d}
\end{figure}

In the second example we regard a pressure wave in straight cylinder as presented in \cite{Gerbeau2003}. The cylinder has the length $\SI{5}{\cm}$ and a radius of $\SI{0.5}{\cm}$. The cylinder is surrounded by an elastic structure with constant thickness $h=\SI{0.1}{\cm}$.
The elastic structure is clamped at the inflow boundary and the outflow domain is fixed in x-direction and free to move in y- and z-direction. At the inlet we describe a pressure step function $p_{in}=1.33\cdot10^4 \SI{}{ \g\per \cm\s}=\SI{10}{\mmHg}$ for $t\leq \SI{0.003}{\s}$, afterwards we set the pressure to zero. In the solid the Lame parameters with $\lambda=1.73\cdot 10^6 \SI{}{\g\per \cm \s^2}$ and $\mu=1.15\cdot 10^6\SI{}{\g\per \cm\s^2}$ and a fluid viscosity $\nu_f=\SI{0.03}{\cm^2\per \s}$ are chosen. The solid and fluid densities are given by $\rho_s=\SI{1.2}{\g \per \cm^3}$ and $\rho_f=\SI{1.0}{\g\per \cm^3}$. We plotted the solution at $t=\SI{0.006}{\s}$ in Figure \ref{Solution_3d}.

If only a do-nothing condition is described with constant pressure at the outflow boundary, then pressure waves are going to be reflected at the outflow boundary. If the pressure along the outflow boundary is chosen appropriate the pressure wave will leave the cylindrical domain without any reflection such that the system will be at rest after some time. To determine the corresponding pressure profile $q(t)\in L^2(I)$ on the time interval $I=[0,0.04]$, we minimize the kinetic energy in the fluid domain for $t>\SI{0.03}{\s}$. Hence we minimize the functional
\begin{equation}
\begin{aligned}
\min_{q \in L^2(I)} J(q,\vt)= \frac{1}{2}\int_{0.03}^{0.04}\norm{\vt}^2_{\FL}\text{ d}t+ \frac{\alpha}{2}\int_{0}^{0.04} q(t)^2 \text{ d}t
\end{aligned}
\end{equation}
constrained by the fluid-structure interaction problem. In time we use, as already in the previous example, a shifted Crank-Nicolson time stepping scheme with time step size $k=0.0001$ and $\theta=0.5+2k$.
Only at the time points $t=0$ and $t=0.003$, we use for four steps a time step-size of $k=0.00005$ with $\theta=1.0$. Thereby, no artificial effects occur due to the jump in the pressure at the inflow boundary.  
The Tikhonov regularization parameter is set to $\alpha=1.0\cdot10^{-8}$.


\begin{table}[t]
  \centering
  \begin{tabular}{l|ccccc}
    \toprule
    mesh level  & 2 & 3 &4 &5 & 6 \\ \midrule
    spatial dofs 2d  	&5\,540	&21\,320	&83\,600
    &331\,040	&1\,317\,440\\
    time steps & \multicolumn{5}{l}{$N=600$ uniform time steps
      (Crank-Nicolson)}\\ 
    \midrule
    spatial dofs 3d  	&43\,904	&336\,224& 894\,656&
    3\,122\,560&-\\
    time steps & \multicolumn{5}{l}{$N=400$ uniform time steps
      (Crank-Nicolson) plus 8 (back. Euler)}\\
    \bottomrule
  \end{tabular}
  \caption{Spatial degrees of freedom for 2d and 3d configuration on
    every refinement level. In 3d the mesh on level 4 and 5 are
    locally refined along the interface. In time we use a uniform
    partitioning. In 3d, we add 4 initial backward Euler steps with
    reduced step size for smoothing at times $t=\unit[0]{s}$ and
    $t=\unit[0.03]{s}$ each. }   
  \label{table:dofs}
\end{table}
%
\subsection{Optimization algorithm}
Given the computed gradient information using Formula
\eqref{eq_gradient}, we apply a BFGS algorithm implemented in the
optimization library RoDoBo \cite{RoDoBo} to solve the optimization
problem. { We use a limited memory version as presented e.g. in \cite{Geiger2013} such that there is no need to
assemble and store the BFGS update matrix. To guarantee that the update matrix keeps symmetric and positive definite a Powell-Wolfe step size control should be used. As this step size criteria is in most cases very conservative, we only check if we have descend in the functional value.} Control constraints could
be added in the optimization algorithm via projection of the update in
every optimization step. In this paper we only regard unconstrained
examples. Fast convergence of the BFGS algorithm only can be expected
close to the optimal solution. Hence, we take advantage of the mesh
hierarchy, which is used in the geometric multigrid algorithm, and
first solve the optimization problem on a coarse grid to have a good
initial guess on finer meshes. As the computation time rises for 3d
examples on finer meshes very fast, we can save a lot of  computation
time using this approach. In the following we reduce the norm of the
gradient by a factor of $10^{-1}$ in every optimization loop and then
refine the mesh and restart the optimization with the control from the
coarser mesh. To compute the gradient, we have to solve the state and
dual problem for all time steps. The state solutions are stored on the hard disc and loaded during the computation of the dual problem if necessary. Thereby, only the current state and adjoint solutions have to be held in the memory.

In every Newton step a GMRES iterative solver preconditioned with a
geometric multigrid solver provided by the FEM software library
Gascoigne \cite{Gascoigne3D} is used. All linear systems are solved up
to a relative accuracy of $10^{-4}$. {In every time step the initial Newton residuum is reduced by a factor of $10^{-6}$.} We use the same tolerances for the state and dual problem. The matrices are only reassembled, if the nonlinear convergence rate falls below $\gamma=0.05$ as in \cite{FailerRichter2019}. In every dual step, the matrices are assembled at least once at the beginning of every Richardson iteration.

\subsection{Numerical results 2d example}
\begin{figure}[t]
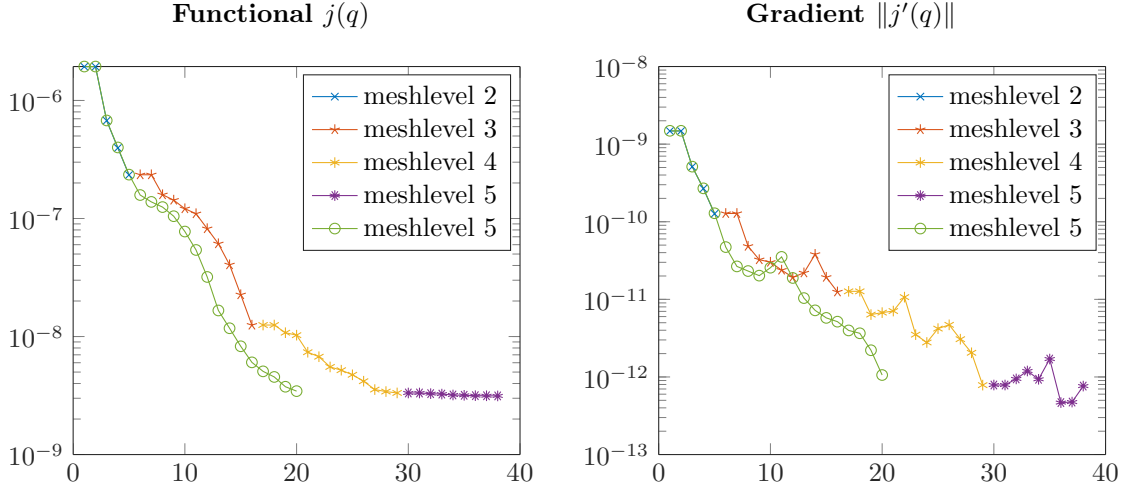

  \centering
  \setlength\figureheight{0.35\textwidth} 
  \setlength\figurewidth{0.40\textwidth} 
  \small 
\begin{tabular}{cc}
    \textbf{Functional $j(q)$}& \textbf{Gradient $\|j'(q)\|$ }\\
\end{tabular}
  \caption{2d example: Functional value plotted over optimization steps (left), Norm of the gradient plotted over optimization steps (right)}
  \label{fig:opt_2d_functional}
\end{figure}

\begin{figure}[t]
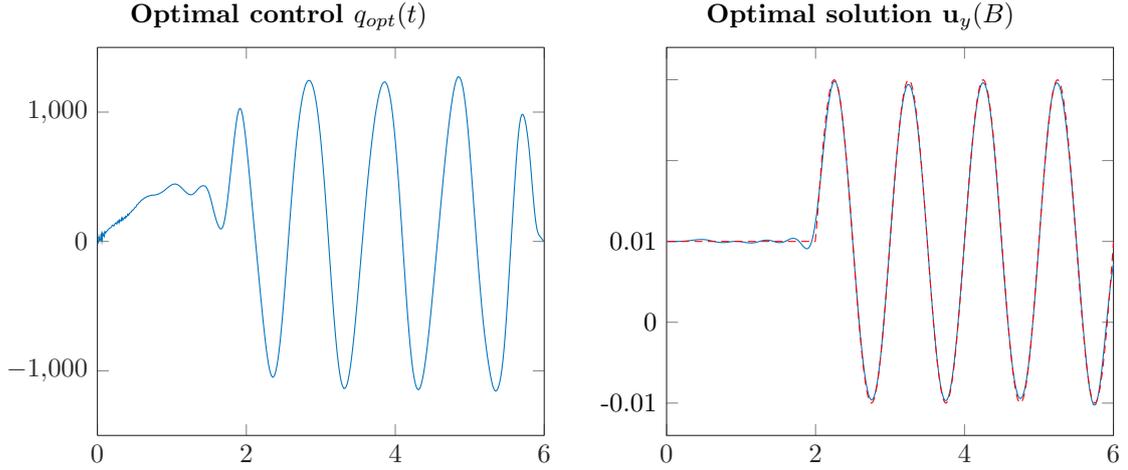

  \centering
  \setlength\figureheight{0.35\textwidth} 
  \setlength\figurewidth{0.40\textwidth} 
  \small 
  \begin{tabular}{cc}
    \textbf{Optimal control $q_{opt}(t)$}& \textbf{Optimal solution
      $\ut_y(B)$ }\\
\end{tabular}
  \caption{2d example: Optimal control $q_{opt}$ plotted over time (left) and optimal solution and reference solution in the point B plotted over time (right)}
  \label{fig:opt_sol_2d}
\end{figure}

In Figure \ref{fig:opt_2d_functional} we plot the value of the regularized functional $j(q)$ and the norm of the gradient $\norm{j'(q)}$ in every optimization step. The optimization algorithm is started with $q_n=0$ for $n=1,\cdots,N$. The computed optimal control is given in Figure~\ref{fig:opt_sol_2d}. The functional value reduces in every optimization step and the gradient can be reduced by a factor of $10^{-3}$ after less then 40 optimization cycles (see Figure~\ref{fig:opt_2d_functional}). In addition we plotted the optimal solution for the 2d example on the finest mesh at the point B in the center of the observation domain $\Omega_{obs}$ and compare the solution to the reference solution in Figure \ref{fig:opt_sol_2d}. Only around the kink of the desired state a mismatch between desired state and optimized solution can be seen.

To evaluate our approach to solve the optimization algorithm first on
coarser grids and then to refine systematically, we restarted the
optimization algorithm directly on meshlevel 5. A similar behavior in
functional and gradient to the previous approach can be observed in
Figure \ref{fig:opt_2d_functional}. To reduce the gradient to a
tolerance of $10^{-12}$ only 20 optimization loops are required. But
since  about $\SI{14000}{\s}$ of computational time are needed to solve one cycle
of the state and the adjoint system for all time steps on meshlevel 5,
but only less then $\SI{4000}{\s}$ on the coarser meshlevel 4  (even
less time on still coarser grids), systematical refinement of the
mesh is much more efficient. The computations were performed on a
Intel(R) Xeon(R) Gold 6150 CPU @ 2.70GHz with 18 threads. We
parallelized the assembling of matrix and the Vanka smoother as well
as the matrix vector multiplication via OpenMP,
see~\cite{FailerRichter2019}.

\subsection{Numerical Results 3d example}

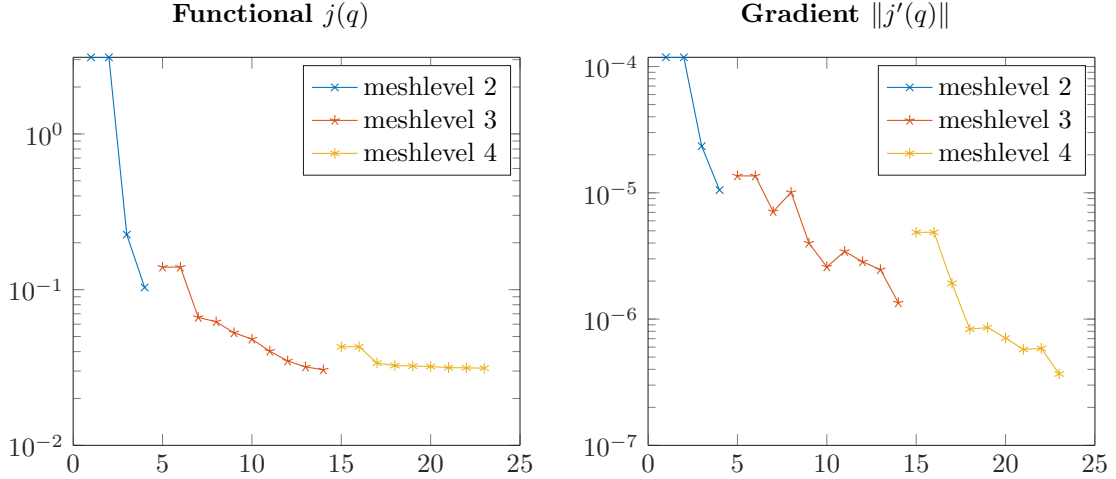
\begin{figure}[t]
  \centering
  \setlength\figureheight{0.35\textwidth} 
  \setlength\figurewidth{0.40\textwidth} 
  \small 
  \begin{tabular}{cc}
   \textbf{Functional $j(q)$}& \textbf{Gradient $\|j'(q)\|$ }\\[2mm]
%
%
%
\definecolor{mycolor1}{rgb}{0,0.447,0.741}%
\definecolor{mycolor2}{rgb}{0.85,0.325,0.098}%
\definecolor{mycolor3}{rgb}{0.929,0.694,0.125}%
\begin{tikzpicture}

\begin{axis}[%
width=\figurewidth,
height=\figureheight,
scale only axis,
separate axis lines,
every outer x axis line/.append style={darkgray!60!black},
every x tick label/.append style={font=\color{darkgray!60!black}},
xmin=0,
xmax=25,
every outer y axis line/.append style={darkgray!60!black},
every y tick label/.append style={font=\color{darkgray!60!black}},
ymode=log,
ymin=0.01,
ymax=3.10336480441781,
yminorticks=true,
legend style={draw=darkgray!60!black,fill=white,legend cell align=left}
]
\addplot [
color=mycolor1,
solid,
mark=x,
mark options={solid}
]
table[row sep=crcr]{
1 3.10336480441781\\
2 3.10336340763207\\
3 0.225580248979941\\
4 0.10321466328132\\
};
\addplot [
color=mycolor2,
solid,
mark=star,
mark options={solid}
]
table[row sep=crcr]{
5 0.139563153471029\\
6 0.139563134846614\\
7 0.0663460161847211\\
8 0.0621739869488137\\
9 0.0528675981923924\\
10 0.0481181835678647\\
11 0.0402534556565917\\
12 0.0348328451119022\\
13 0.0318980818276145\\
14 0.0305731735226139\\
};
\addplot [
color=mycolor3,
solid,
mark=asterisk,
mark options={solid}
]
table[row sep=crcr]{
15 0.042977738131788\\
16 0.0429777356914979\\
17 0.0337566850098577\\
18 0.0325351879239444\\
19 0.0323171780279687\\
20 0.0320914217938281\\
21 0.0315691548172662\\
22 0.031433066922041\\
23 0.0313276959717998\\
};
\addlegendentry{meshlevel 2}
\addlegendentry{meshlevel 3}
\addlegendentry{meshlevel 4}
\end{axis}
\end{tikzpicture}%

   &
%
%
%
\definecolor{mycolor1}{rgb}{0,0.447,0.741}%
\definecolor{mycolor2}{rgb}{0.85,0.325,0.098}%
\definecolor{mycolor3}{rgb}{0.929,0.694,0.125}%
\begin{tikzpicture}

\begin{axis}[%
width=\figurewidth,
height=\figureheight,
scale only axis,
separate axis lines,
every outer x axis line/.append style={darkgray!60!black},
every x tick label/.append style={font=\color{darkgray!60!black}},
xmin=0,
xmax=25,
every outer y axis line/.append style={darkgray!60!black},
every y tick label/.append style={font=\color{darkgray!60!black}},
ymode=log,
ymin=1e-07,
ymax=0.00011826180707105,
yminorticks=true,
legend style={draw=darkgray!60!black,fill=white,legend cell align=left}
]
\addplot [
color=mycolor1,
solid,
mark=x,
mark options={solid}
]
table[row sep=crcr]{
1 0.00011826180707105\\
2 0.000118261771235517\\
3 2.34115703915687e-05\\
4 1.05128222529164e-05\\
};
\addplot [
color=mycolor2,
solid,
mark=star,
mark options={solid}
]
table[row sep=crcr]{
5 1.36118228904098e-05\\
6 1.36118208103327e-05\\
7 7.10473954372292e-06\\
8 1.00784644472018e-05\\
9 3.99995461506839e-06\\
10 2.58985009923667e-06\\
11 3.4370615125192e-06\\
12 2.84471212117245e-06\\
13 2.46352857052686e-06\\
14 1.34633223753749e-06\\
};
\addplot [
color=mycolor3,
solid,
mark=asterisk,
mark options={solid}
]
table[row sep=crcr]{
15 4.86539714181047e-06\\
16 4.86539645694548e-06\\
17 1.91857760814547e-06\\
18 8.31893855782116e-07\\
19 8.57094441595686e-07\\
20 7.07868982982835e-07\\
21 5.74854456716177e-07\\
22 5.86320384694015e-07\\
23 3.68521476618369e-07\\
};
\addlegendentry{meshlevel 2}
\addlegendentry{meshlevel 3}
\addlegendentry{meshlevel 4}
\end{axis}
\end{tikzpicture}%

  \end{tabular}
  \caption{3d example: Functional value plotted over optimization steps  (left), Norm of the gradient plotted over optimization steps (right)}
  \label{fig:opt_3d_functional}
\end{figure}

For the 3d example we used the same optimization algorithm as in the 2d case. We can see in Figure \ref{fig:opt_3d_functional} that using the optimized pressure profile at the outflow boundary about $98.9\%$ of the kinetic energy after $t>0.03s$ now leaves the cylinder. The jumps in the gradient after every refinement step indicate that a more accurate computation on the coarse grid would not result in better starting values on the finer meshes. The norm of the gradient could be reduced from $1.15\cdot 10^{-4}$ to $3.68\cdot 10^{-7}$ in 23 optimization steps, whereby only 9 optimization cycles were necessary on the computationally costly fine grid on meshlevel 4.

\begin{figure}[t]
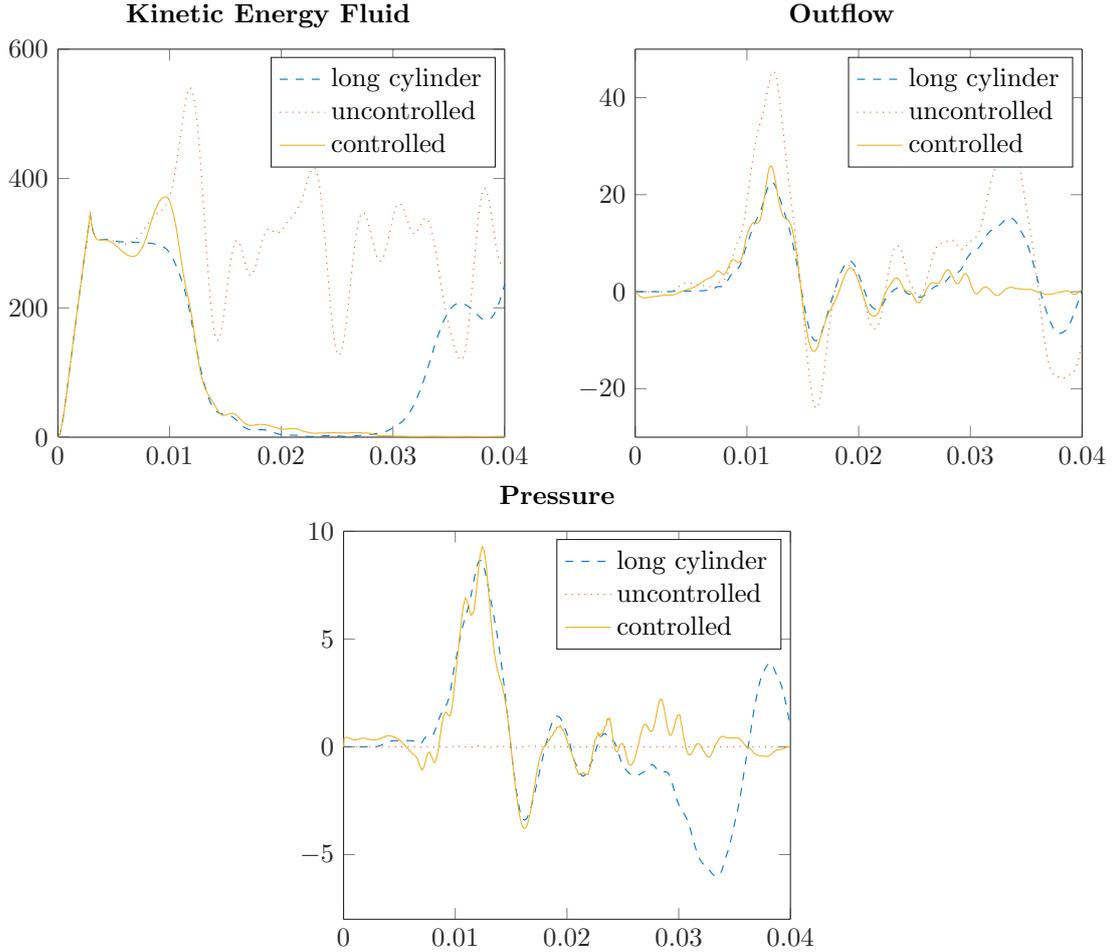

  \centering
  \setlength\figureheight{0.35\textwidth} 
  \setlength\figurewidth{0.40\textwidth} 
  \small
\begin{tabular}{cc}
    \textbf{Kinetic Energy Fluid}& \textbf{Outflow }\\
\end{tabular}
  \caption{Kinetic energy in the fluid domain $\FL$ as well as outflow and mean pressure plotted over time at $\Gamma_q$ for $q=0$, $q_{opt}$ and for a long cylinder}
  \label{fig:opt3d_Energy}
\end{figure}

To compare the controlled solution with the uncontrolled solution, we computed in addition the solution on a cylinder with length $\SI{10}{\cm}$. As the reflection on the outflow boundary occurs later the pressure and flow values in the center at $x=\SI{5}{\cm}$ can be seen as reference values for optimal non reflective boundary conditions for $t<\SI{0.02}{\s}$. As we only control the pressure on the fluid domain, reflections on the solid boundary can still occur. Furthermore we can observe that the pressure is not constant along the virtual outflow boundary at $x=\SI{5}{\cm}$ for the long cylinder. Thus, we can not expect the reference solution to fully match the optimized solution. Nevertheless, we can see in Figure \ref{fig:opt3d_Energy} that pressure and outflow profiles of the controlled solution are very close to the reference values at the outflow boundary. In addition the kinetic energy in the fluid domain has a similar decay behavior. After the time point $t=\SI{0.025}{\s}$ the kinetic energy in the left half of the long cylinder rises again due to the reflection of the pressure wave at the outflow boundary at $x=\SI{10}{\cm}$. This explains the different behavior of pressure and outflow after $t=\SI{0.025}{\s}$.

\subsection{Test of the Newton scheme, of the Richardson iteration and of the iterative linear solver}
How to evaluate the performance of the quasi Newton scheme or of the
iterative linear solver is not so obvious. Due to the changing control
in every optimization cycle and the nonlinear character of the
problem, the condition numbers of the matrices occurring in the linear
subproblems will vary in every time step and optimization
cycle. Hence, we first compute only one optimization step with $q=0$
on various meshlevels to analyze the h-dependence of our solution
algorithm. Thereby, we compare the mean number of Newton/Richardson
iterations and GMRES steps per time step on every meshlevel. Afterwards, we compute mean values in every optimization loop to analyze how the performance can vary during the optimization procedure. 

\begin{table}[t]
  \centering
%
  \begin{tabular}{c|cc|ccc}    
    \toprule
    mesh & Newton- & Matrix- &
    GMRES~(\ref{problem:1}) & GMRES~(\ref{problem:2}) &
    GMRES~(\ref{problem:3}) \\
    level & steps & assemblies & (momentum) & (deformation) & (extension) \\
    \midrule
    3&3.75&0.28&7.19&1.26&3.77\\
    4&3.60&0.87&8.34&1.27&3.82\\
    5&3.91&1.07&9.52&1.25&3.77\\
    6&4.33&1.54&10.61&1.31&3.98\\
    \bottomrule
    \multicolumn{6}{c}{~}\\ 
    \toprule
    mesh & Richardson- & Matrix- &
    GMRES~(\ref{equ:dualproblem1}) & GMRES~(\ref{equ:dualproblem2}) &
    GMRES~(\ref{equ:dualproblem3}) \\
    level & steps & assemblies & (extension) & (deformation) & (momentum) \\
    \midrule
    3&3.09&1.20&4.38&1.32&7.88\\
    4&3.50&1.43&4.80&1.28&8.68\\
    5&3.64&1.48&5.86&1.27&9.54\\
    6&3.69&1.52&6.63&1.27&11.17\\
    \bottomrule
\end{tabular}
  \caption{Average number of Newton/Richardson iterations, average
    number of matrix assemblies per time
    step and average number of GMRES steps for solving the three
    subproblems in the first optimization step. Top: 2d example
    state,~(\ref{problem:1}) is the coupled momentum
    equation,~(\ref{problem:2}) the coupling between  solid velocity
    and deformation and~(\ref{problem:3}) the fluid 
    deformation extension. Bottom: 2d
    example dual, where~(\ref{equ:dualproblem1}) is the adjoint
    extension equation,~(\ref{equ:dualproblem2}) the adjoint solid
    velocity-deformation coupling and~(\ref{equ:dualproblem3}) the
    adjoint coupled momentum equation.}
  \label{table:Newtonsteps2dTime}
\end{table}

In Table \ref{table:Newtonsteps2dTime} and  Table
\ref{table:Newtonsteps3dTime} we present the mean number of
Newton/Richardson iterations and Matrix assemblies per time step for the 2d and 3d
examples. In addition we present the mean  number of GMRES steps needed
per Newton/Richardson iteration to solve the linear subproblems \eqref{problem:3}, \eqref{problem:2} and \eqref{problem:1} and 
\eqref{equ:dualproblem1}, \eqref{equ:dualproblem2} and
\eqref{equ:dualproblem3}. We can observe that the number of
Newton/Richardson iterations per time step ranges between 3 and 4 for
state and dual problem. The coupled momentum equation remains the most
complex problem with the highest number of steps
required. Equations~(\ref{problem:2}) and~(\ref{equ:dualproblem2})
belong to the discretization of the velocity deformation coupling $d_t
u=v$ within the solid domain. This corresponds to the inversion of the
mass matrix which explains the very low iteration counts. 
The results for the state equation are similar
to the values in \cite{FailerRichter2019}, where we already could
observe for different examples that neglecting the ALE derivatives
only has minor influences on the behavior of the Newton scheme. In
\cite{FailerRichter2019} a more detailed analysis of the smoother in
the geometric multigrid algorithm can be found.

As the dual equation
is linear with respect to the adjoint variable, we would have expected
to need only one Richardson iteration per time step. However, since we
neglected terms occurring due to the ALE transformation, we loose the
optimal convergence and need about 3 Richardson iterations per
time step. On the other hand, the matrices occurring  in cascade of
subproblems in the dual system have the same condition number as the
matrices in the state equation. Only by this approximation and
splitting, iterative solvers can successfully be applied to solve the
linear problems. The number of GMRES steps in the dual subproblems in
Table \ref{table:Newtonsteps2dTime} and  Table
\ref{table:Newtonsteps3dTime} are rather small and close to the values
for the state problem. As \eqref{equ:dualproblem3} and
\eqref{problem:1} are still fully coupled problems of fluid and
elastic solid, most of the computational time is spent in solving
these two subproblems. The least number of GMRES steps is needed to
invert the mass matrix in \eqref{problem:2} and
\eqref{equ:dualproblem2}. In all subproblems the number of GMRES steps
only increases slightly under mesh refinement. 

\begin{figure}[t]
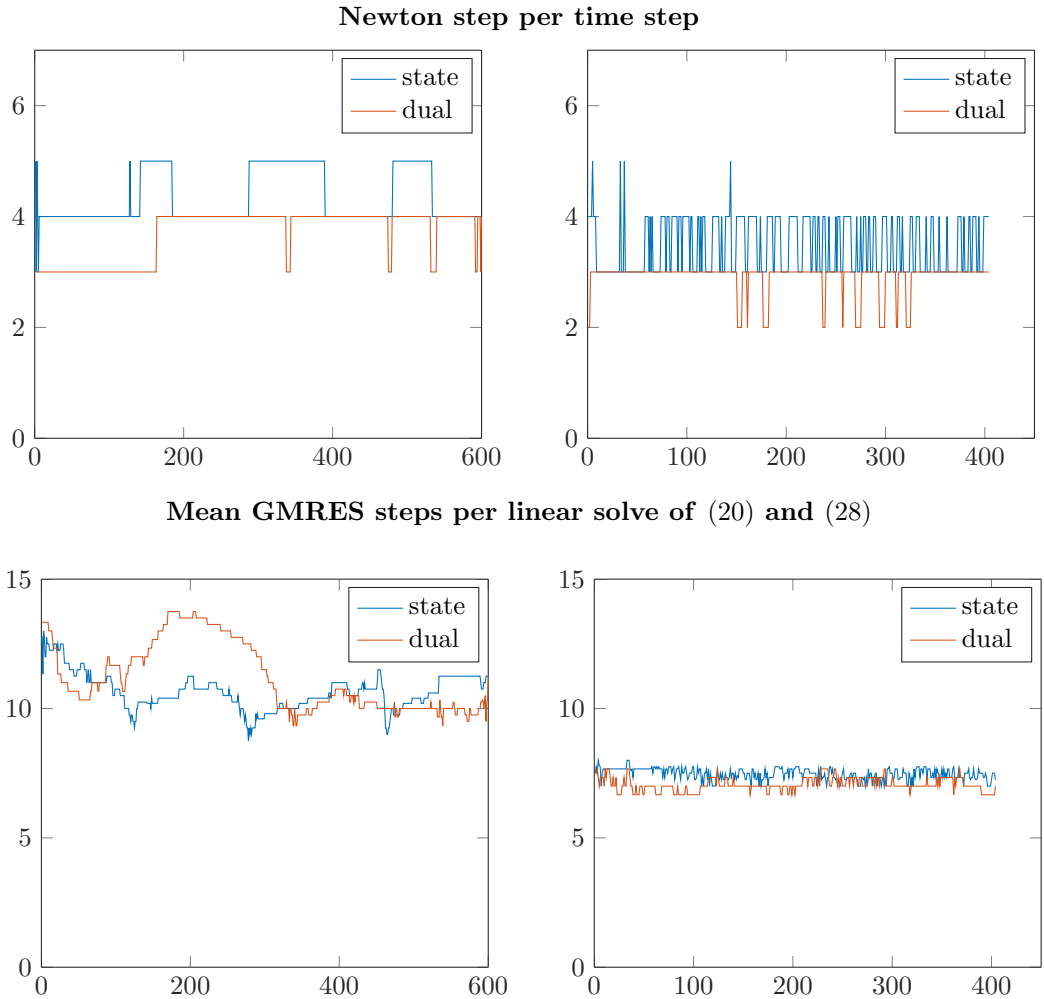

  \centering
  \setlength\figureheight{0.35\textwidth} 
  \setlength\figurewidth{0.40\textwidth} 
  \small
\begin{tabular}{cc}
    \multicolumn{2}{c}{\textbf{Newton step per time step}}\\
\end{tabular}
  \caption{Performance in first optimization loop on mesh level 5 in 2d (Left) and on mesh level 4 in 3d (Right).
   Mean GMRES steps per linear solve of \eqref{problem:1} and
   \eqref{equ:dualproblem3} plotted over time steps and number of
   Newton steps plotted over time steps}
  \label{fig:StateDualPerform2d}
\end{figure}

In Figure \ref{fig:StateDualPerform2d} we show the mean number of
GMRES steps per Newton step in every timestep and the number of Newton
steps per timestep. For the given examples the values only vary
slightly over time. Therefore the mean values in Table
\ref{table:Newtonsteps2dTime} and  Table \ref{table:Newtonsteps3dTime}
represent the overall behavior very well.

\begin{figure}[t]
  \centering
  \setlength\figureheight{0.35\textwidth} 
  \setlength\figurewidth{0.40\textwidth} 
  \small 
%
%
%
\definecolor{mycolor1}{rgb}{0,0.447,0.741}%
\definecolor{mycolor2}{rgb}{0.85,0.325,0.098}%
\begin{tikzpicture}

\begin{axis}[%
width=\figurewidth,
height=\figureheight,
scale only axis,
separate axis lines,
every outer x axis line/.append style={darkgray!60!black},
every x tick label/.append style={font=\color{darkgray!60!black}},
xmin=0,
xmax=20,
every outer y axis line/.append style={darkgray!60!black},
every y tick label/.append style={font=\color{darkgray!60!black}},
ymin=6,
ymax=12,
legend style={draw=darkgray!60!black,fill=white,legend cell align=left}
]
\addplot [
color=mycolor1,
solid,
mark=o,
mark options={solid}
]
table[row sep=crcr]{
1 9.37653018151119\\
2 9.37441958632334\\
3 9.28971553610503\\
4 9.32471008028546\\
5 9.32266549758878\\
6 9.35235624729788\\
7 9.28503460207613\\
8 9.21061792863359\\
9 9.10331632653061\\
10 8.94991582491583\\
11 8.77992440151197\\
12 8.75073313782991\\
13 8.70600500417014\\
14 8.7029497299543\\
15 8.65769867549669\\
16 8.64173553719008\\
17 8.63331955353452\\
18 8.61207609594706\\
19 8.62950683796104\\
20 8.63098708871304\\
};
\addplot [
color=mycolor2,
solid,
mark=+,
mark options={solid}
]
table[row sep=crcr]{
1 9.27885921231326\\
2 9.27885921231326\\
3 8.84187082405345\\
4 8.69779735682819\\
5 8.5792656587473\\
6 8.41858482523444\\
7 8.41354344122658\\
8 8.44776119402985\\
9 8.57480980557904\\
10 8.52866779089376\\
11 8.54153225806452\\
12 8.57736905237905\\
13 8.5622283682339\\
14 8.59213089209194\\
15 8.62681871804955\\
16 8.64656549520767\\
17 8.62475364603863\\
18 8.63130711164084\\
19 8.64915719325755\\
20 8.64274447949527\\
};
\legend{state,dual}
\end{axis}
\end{tikzpicture}%
%
%
%
%
\definecolor{mycolor1}{rgb}{0,0.447,0.741}%
\definecolor{mycolor2}{rgb}{0.85,0.325,0.098}%
\begin{tikzpicture}

\begin{axis}[%
width=\figurewidth,
height=\figureheight,
scale only axis,
separate axis lines,
every outer x axis line/.append style={darkgray!60!black},
every x tick label/.append style={font=\color{darkgray!60!black}},
xmin=0,
xmax=20,
every outer y axis line/.append style={darkgray!60!black},
every y tick label/.append style={font=\color{darkgray!60!black}},
ymin=3,
ymax=5,
legend style={draw=darkgray!60!black,fill=white,legend cell align=left}
]
\addplot [
color=mycolor1,
solid,
mark=o,
mark options={solid}
]
table[row sep=crcr]{
1 3.94833333333333\\
2 3.94833333333333\\
3 3.80833333333333\\
4 3.73666666666667\\
5 3.80166666666667\\
6 3.855\\
7 3.85333333333333\\
8 3.83\\
9 3.92\\
10 3.96\\
11 3.96833333333333\\
12 3.97833333333333\\
13 3.99666666666667\\
14 4.01166666666667\\
15 4.02666666666667\\
16 4.03333333333333\\
17 4.03166666666667\\
18 4.03\\
19 4.02166666666667\\
20 4.00166666666667\\
};
\addplot [
color=mycolor2,
solid,
mark=+,
mark options={solid}
]
table[row sep=crcr]{
1 3.68166666666667\\
2 3.68166666666667\\
3 3.74166666666667\\
4 3.78333333333333\\
5 3.85833333333333\\
6 3.91\\
7 3.91333333333333\\
8 3.90833333333333\\
9 3.94333333333333\\
10 3.95333333333333\\
11 4.13333333333333\\
12 4.16833333333333\\
13 4.21833333333333\\
14 4.27833333333333\\
15 4.23833333333333\\
16 4.17333333333333\\
17 4.22833333333333\\
18 4.195\\
19 4.25166666666667\\
20 4.22666666666667\\
};
\legend{state,dual}
\end{axis}
\end{tikzpicture}%
  \caption{Performance in every optimization loop on mesh level 5 in
    2d. Left: Mean GMRES steps per linear solve of \eqref{problem:1}
    and \eqref{equ:dualproblem3} plotted over optimization steps
    . Right: Mean Newton steps per time step plotted over optimization steps}
  \label{fig:opt_2d_gmres_per_optloop}
\end{figure}
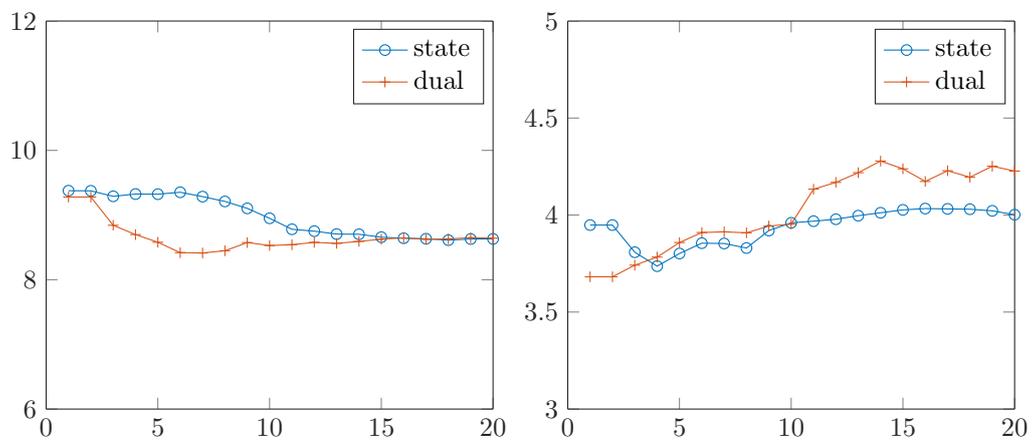

To understand how the performance of the solution algorithm in the
case of the 2d example varies during the optimization loop, we
show  in Figure \ref{fig:opt_2d_gmres_per_optloop} the average number
of Newton steps per time step and the average number of GMRES steps
per Newton step for each optimization step. The computation was
started with $q=0$ on meshlevel 5. No further mesh was applied.  

The dependency on the time step size of the presented quasi Newton
scheme for the state equation was already analyzed in
\cite{FailerRichter2019}. Therein, we could observe an increasing
Newton iteration count for larger time steps, in particular for very
large time steps. A similar behavior can be expected for the dual
problem. While the presented solution approach can be regarded as
highly efficient and appropriate for optimal control of nonstationary
fluid-structure interactions, stationary or quasi stationary
applications will call for alternative approaches like the geometric
multigrid solver presented by Aulisa, Bna and Bornia~\cite{AulisaBnaBornia2018}.

In \cite{FailerRichter2019} more details regarding the computational
time and savings due to parallelization can be found. As we have to
evaluate state and adjoint variables, as well as additional terms due
to linearization in every dual step, the computational time to
assemble the matrix and to compute the Newton residuum is slightly
larger for the dual equation then for the state equation.

\begin{table}[t]
  \centering
%
  \begin{tabular}{c|cc|ccc}    
    \toprule
    mesh & Newton- & Matrix- &
    GMRES~(\ref{problem:1}) & GMRES~(\ref{problem:2}) &
    GMRES~(\ref{problem:3}) \\
    level & steps & assemblies & (momentum) & (deformation) & (extension) \\
    \midrule
    2&4.21&0.09&6.82&1.46&2.47\\
    3&3.35&0.80&7.04&1.95&2.60\\
    4&3.46&0.52&7.47&1.91&2.57\\
    5&3.70&0.48&7.34&1.83&2.55\\
    \bottomrule
    \multicolumn{6}{c}{~}\\ 
    \toprule
    mesh & Richardson- & Matrix- &
    GMRES~(\ref{equ:dualproblem1}) & GMRES~(\ref{equ:dualproblem2}) &
    GMRES~(\ref{equ:dualproblem3}) \\
    level & steps & assemblies & (extension) & (deformation) & (momentum) \\
    \midrule
    2&2.99&1.00&2.71&1.35&6.74\\
    3&3.00&1.00&3.30&1.59&7.03\\
    4&2.90&1.00&3.32&1.47&7.06\\
    5&2.39&1.00&3.63&1.51&7.01\\
    \bottomrule
\end{tabular}
  \caption{Average number of Newton/Richardson iterations, average
    number of matrix assemblies per time
    step and average number of GMRES steps for solving the three
    subproblems in the first optimization step. Top: 3d example
    state,~(\ref{problem:1}) is the coupled momentum
    equation,~(\ref{problem:2}) the coupling between  solid velocity
    and deformation and~(\ref{problem:3}) the fluid 
    deformation extension. Bottom: 3d
    example dual, where~(\ref{equ:dualproblem1}) is the adjoint
    extension equation,~(\ref{equ:dualproblem2}) the adjoint solid
    velocity-deformation coupling and~(\ref{equ:dualproblem3}) the
    adjoint coupled momentum equation.}
  \label{table:Newtonsteps3dTime}
\end{table}

\section{Summary}\label{sec:conclusion}

We extended the Newton multigrid framework for monolithic
fluid-structure interactions in ALE coordinates presented in
\cite{FailerRichter2019} to the dual system of fluid-structure
interaction problems. The solver is based on replacing the adjoint
solid deformation by a new variable and on skipping the ALE
derivatives within the adjoint Navier-Stokes equations. As we do not
modify the residuum, state and dual solution in
each time step still converge to the exact discrete solution.  The
adjoint coupling conditions incorporated in the monolithic formulation
are still fulfilled. As we compute correct gradient information, we
see fast convergence in our optimization algorithm. The coupled system
is better conditioned (as compared to monolithic Jacobians) which
allows to use very simple multigrid smoothers that are easy to
parallelize. Only this makes gradient based algorithm feasible and
efficient for 3d fluid-structure interaction applications, where memory consumption prevents the use of direct solvers.

It would be straightforward to combine the presented algorithm with
dual-weighted residual error estimators for mesh and time step refinement. Instead of global refinement of the mesh after every optimization loop the error estimators indicate where to refine the mesh locally. The sensitivity information from the optimization algorithm can be directly used to evaluate the error estimators.
 
%
%

\section*{Acknowledgements}
Both authors acknowledge the financial support by the Federal Ministry
of Education and Research of Germany, grant number 05M16NMA (TR) and 05M16WOA
(LF). TR further acknowledges the support of the GRK 2297 MathCoRe, funded by the
Deutsche Forschungsgemeinschaft, grant number 314838170 and LF acknowledges the support
of the GRK 1754, funded by the Deutsche Forschungsgemeinschaft, grant number 188264188.

\end{document}